\documentclass[11pt]{article}
\usepackage{amsfonts}
\usepackage{mathrsfs}
\usepackage{cite}
\usepackage{amssymb,amsmath} % font used for R in Real numbers
 \allowdisplaybreaks

\catcode`!=11
\let\!int\int \def\int{\displaystyle\!int}
\let\!lim\lim \def\lim{\displaystyle\!lim}
\let\!sum\sum \def\sum{\displaystyle\!sum}
\let\!sup\sup \def\sup{\displaystyle\!sup}
\let\!inf\inf \def\inf{\displaystyle\!inf}
\let\!cap\cap \def\cap{\displaystyle\!cap}
\let\!max\max \def\max{\displaystyle\!max}
\let\!min\min \def\min{\displaystyle\!min}
\let\!frac\frac \def\frac{\displaystyle\!frac}
\catcode`!=12

\let\oldsection\section
\renewcommand\section{\setcounter{equation}{0}\oldsection}

\allowdisplaybreaks
\def\pf{\it{Proof.}\rm\quad}

\def\R{\mathbb{R}}

\newtheorem{thm}{Theorem}[section]

\newtheorem{pro}{Proposition}[section]
\newtheorem{lem}{Lemma}[section]
\newtheorem{re}{Remark}[section]
\begin{document}
%%%%%%%%%%%%%%%%%%%% title %%%%%%%%%%%%%%%%%%%%%%%%%%%%%%%%%%%%%%%%%%%%%%%%
\title{\Large\bf THE STABILITY OF STRONG VISCOUS  CONTACT DISCONTINUITY TO AN INFLOW PROBLEM FOR FULL COMPRESSIBLE NAVIER-STOKES EQUATIONS     }
\author{{Tingting Zheng\footnote{Corresponding author, e-mail: asting16@sohu.com}}\\[2mm]
{\small  Computer and Message Science College,
 Fujian Agriculture
and Forestry University,}\\[2mm]
{\small  Fuzhou 350001, P. R. China }}

\date{}

\maketitle

\noindent{\bf Abstract.} This paper is concerned with nonlinear
stability of viscous contact discontinuity to inflow problem for the
one-dimensional full compressible Navier-Stokes equations with
different ends in half space $[0,\infty)$. For the case when the
local stability of the contact discontinuities was first studied by
\cite{X},later generalized by \cite{LX}, local stability of weak
viscous contact discontinuity is well-established by
\cite{HMS,HMX,HXY,HZ,HLM2009}, but for the global stability of
inflow gas with big oscillation ends $(|\theta_+-\theta_-|>1\ and \
|\rho_+-\rho_-|>1)$, fewer results have been obtained excluding zero
dissipation \cite{MaSX} or $\gamma\to 1$ gas see \cite{HH}. Our main
purpose is to deduce the corresponding nonlinear stability result
with the two different ends by exploiting the elementary energy
method. As a first step towards this goal, we will show in this
paper that with a certain class of big perturbation which can allow
$|\theta_--\theta_+|>1$ and $|\rho_--\rho_+|>1$ ,the global
stability result holds.
\\[2mm]
\\[2mm]
\section{Introduction}
This paper is concerned with an ``inflow problem" for a
one-dimensional compressible viscous heat-conducting flow in the
half space $\R_+=[0,\infty)$, which is governed by the following
initial-boundary value problem in Eulerian coordinate
$(\tilde{x},t)$:
\begin{equation}
\left\{
\begin{array}{lll}
\tilde\rho_t+(\tilde\rho \tilde u)_{\tilde x}=0,\quad (\tilde x,t)\in\R_+\times\R_+,\\[2mm]
(\tilde\rho \tilde u)_t+(\tilde\rho \tilde u^2+\tilde p)_{\tilde{x}}=\mu \tilde u_{\tilde{x}\tilde{x}}, \\[2mm]
\left(\tilde\rho\left(\tilde e+\frac{\tilde
u^2}{2}\right)\right)_t+\left(\tilde\rho \tilde u\left(\tilde
e+\frac{\tilde u^2}{2}\right)+\tilde p\tilde
u\right)_{\tilde{x}}=\kappa\tilde
\theta_{\tilde{x}\tilde{x}}+(\mu\tilde u\tilde
u_{\tilde{x}})_{\tilde{x}},
\\[4mm]
(\tilde\rho,\tilde u,\tilde \theta)|_{\tilde{x}=0}=(\rho_-,u_b,\theta_-)\quad\mbox{with}\quad u_b>0, \\[2mm]
(\tilde \rho,\tilde u,\tilde \theta)|_{t=0}=(\tilde\rho_0,\tilde
u_0,\tilde\theta_0)(\tilde x)\to(\rho_+,u_b,\theta_+)
\quad\mbox{as}\quad \tilde{x}\to\infty,
\end{array}\right.
\label{1.1}
\end{equation}
where $\tilde\rho$, $\tilde u$ and $\tilde\theta$ are the density,
the velocity and the absolute temperature, respectively, while
$\mu>0$ is the viscosity coefficient and $\kappa>0$ is the
heat-conductivity coefficients, respectively. It is assumed
throughout the paper that $\rho_\pm$, $u_b$ and $\theta_\pm$ are
prescribed positive constants . We shall focus our interest on the
 polytropic ideal gas  with
$|\theta_+-\theta_-|$ and  $|\rho_+-\rho_-|$ are general constants
(not small)    , so  the pressure $\tilde p=\tilde
p(\tilde\rho,\tilde\theta)$ and the internal energy $\tilde e=\tilde
e(\tilde\rho,\tilde\theta)$ are related by the second law of
thermodynamics:
\begin{equation}\label{1.2}
\tilde p=R\tilde\rho\tilde\theta,\quad \tilde
e=\frac{R}{\gamma-1}\tilde\theta+const.,
\end{equation}
where $\gamma>1$ is the adiabatic exponent and $R>0$ is the gas
constant.

The boundary condition $(\ref{1.1})_4$ implies that, through the
 boundary $\tilde
x=0$, the fluid with density $\rho_-$ flows into the region $\R_+$
at the speed $u_b>0$. So the initial-boundary value problem
(\ref{1.1}) is the so-called {\it inflow} problem. On the other
hand, in the case that $u_b=0$ (resp. $u_b<0$), the problem is
called the {\it impermeable wall} (resp. {\it outflow}) problem in
which the boundary condition of density can't be imposed. In terms
of various boundary values,  Matsumura \cite{Ma2001} classified all
possible large-time behaviors of the solutions for the
one-dimensional (isentropic)compressible Navier-Stokes equations.

 Our main purpose is to study the asymptotic stability of the  contact discontinuity
for the inflow problem (\ref{1.1}). It is well known that there are
three basic wave patterns for the $1D$ compressible Euler equations,
including two nonlinear waves,say shock and rarefaction waves, and a
linearly degenerate wave, say contact discontinuity. There have been
a lot of works on the asymptotic behaviors of solutions to the
initial-boundary value (or Cauchy) problem for the Navier-Stokes
equations toward these basic waves or their viscous versions, see,
for example,[3--26]and the reference therein. In what follows, we
briefly recall some related references. Concerning the inflow
problem, Matsumura and Nishihara \cite{MN2001} considered an inflow
problem for the one-dimensional isentropic model system of
compressible viscous gas (i.e.the $1D$ isentropic Navier-Stokes
$(\ref{1.1})_1-(\ref{1.1})_2$ with
$\tilde{p}=R\tilde{\rho}^{\gamma}$) and established the stability
theorems on both the boundary layer solution and the superposition
of a boundary-layer solution and a rarefaction wave. We also refer
to the paper due to Huang et al. \cite{HMS2003} in which the
asymptotic stability on both the viscous shock wave and the
superposition of a viscous shock wave and a boundary-layer solution
are studied. On the other hand, the problem of stability of contact
discontinuities are associated with linear degenerate fields and are
less stable than the nonlinear waves for the inviscid system (Euler
equations). It was observed in \cite{X,LX}, where the metastability
of contact waves was studied for viscous conservation laws with
artificial viscosity, that the contact discontinuity cannot be the
asymptotic state for the viscous system, and a diffusive wave, which
approximated the contact discontinuity on any finite time interval,
actually dominates the large-time behavior of solutions. The
nonlinear stability of contact discontinuity for the (full)
compressible Navier-Stokes equations was then investigated in
\cite{HMS,HZ}
 for the free  boundary value problem and \cite{HMX,HXY}for the Cauchy problem.

However, to our best knowledge, fewer mathematical literature known
for the large-time behaviors of solutions to the inflow problem of
the full compressible viscous heat-conducting Navier-Stokes
equations due to various difficulties come from the big oscillation
ends. So the aim of this paper is to show that the  contact
discontinuities are metastability wave patterns for the inflow
problem (\ref{1.1}) of the full Navier-Stokes system.

To state our main results we first transfer (\ref{1.1}) to the
problem in the Lagrangian coordinate and then make use of a
coordinate transformation to reduce the initial-boundary value
problem (\ref{1.1}) into the following form:
\begin{equation}
\left\{
\begin{array}{ll}
v_t-sv_x-u_x=0,\quad (x,t)\in\R_+\times\R_+,\\[2mm]
u_t-su_x+\left(\frac{R\theta}{v}\right)_x=\mu\left(\frac{u_x}{v}\right)_x,
\\[4mm]
\frac{R}{\gamma-1}\theta_t-\frac{R}{\gamma-1}
s\theta_x+R\frac{\theta}{v}u_x=\kappa\left(\frac{\theta_x}{v}\right)_x+\mu\frac{u_x^2}{v},\\[4mm]
(v,u,\theta)|_{x=0}=(v_-,u_b,\theta_-),\quad t>0,\\[2mm]
(v,u,\theta)|_{t=0}=(v_0,u_0,\theta_0)\to(v_+,u_b,\theta_+)\quad\mbox{as}\quad
x\to\infty,
\end{array}
\right.\label{1.3}
\end{equation}
where $v_\pm$, $u_b$ and $\theta_\pm$ are given positive constants,
and $s=-u_b/v_-<0$, $v_0,\ \theta_0>0$. In fact $v=1/\rho(x,t),\
u=u(x,t),\ \theta=\theta(x,t)$ and $R\theta/v=p(v,\theta)$ are the
specific volume, velocity , temperature and pressure as in
(\ref{1.1}).

Recently most of these Navier-Stokes equations use their Euler
systems as their limitations. Here we consider the
 corresponding Euler system of (\ref{1.3}) with Riemann initial
data reads as follows:
\begin{equation}
\left\{
\begin{array}{ll}
v_t-sv_x-u_x=0,\\[2mm]
u_t-su_x+p(v,\theta)_x=0,\\[2mm]
\frac{R}{\gamma-1}\theta_t-\frac{R}{\gamma-1}
s\theta_x+R\frac{\theta}{v}u_x=0,\\[4mm]
(v,u,\theta)(x,0)=(v_-,u_b,\theta_-)\quad\mbox{if}\quad x<0,\\[2mm]
(v,u,\theta)(x,0)=(v_+,u_b,\theta_+)\quad\mbox{if}\quad x>0,
\end{array}
\right.\label{1.4}
\end{equation}
where $v\pm=\frac{1}{\rho_{\pm}}$, $u_b$ and $\theta\pm$ are the
same positive constants as in (\ref{1.1}).

 Because the corresponding Euler
equations (\ref{1.4}) with the Riemann initial data has the
following soluitons
\begin{equation}
\left(\overline{V},\overline{U},\overline{\Theta}\right) =\left\{
\begin{array}{ll}
(v_-,u_b,\theta_-),\quad x<-st,\\[2mm]
(v_+,u_b,\theta_+),\quad x>-st,
\end{array}
\right. \label{1.5}
\end{equation}
provided that
\begin{equation}
p_-=R\frac{\theta_-}{v_-}=p_+=R\frac{\theta_+}{v_+}.\label{1.6}
\end{equation}
As that in \cite{HMS} we conjecture that the asymptotic limit
$(V,U,\Theta)$ of (\ref{1.3}) is as follows

\begin{equation}
P(V,\Theta)=R\frac{\Theta}{V}=p_+,\;\;\;U(x,t)=\frac{\kappa(\gamma-1)\Theta_x}{\gamma
R\Theta}+u_b, \label{1.7}
\end{equation}
and $\Theta $ is the solution of the following  problem
\begin{equation}
\left\{ \begin{array}{ll}
\Theta_t-s\Theta_x=a(\ln\Theta)_{xx},\quad a=\frac{\kappa p_+(\gamma-1)}{\gamma R^2}>0,\\[2mm]
\Theta(0,t)=\theta_-,\\[2mm]
 \Theta(x,0)=\Theta_{0}\to \theta_+.
\end{array}\right.\label{1.8}
\end{equation}
$(V,U,\Theta)$ satisfies
\begin{equation}
\left\{
\begin{array}{ll}
R\frac{\Theta}{V}=p_+,\\[4mm]
V_t-sV_x=U_x,\\[2mm]
U_t-sU_x+P(V,\Theta)_x=\mu\left(\frac{U_x}{V}\right)_x+F,\\[4mm]
\frac{R}{\gamma-1}\Theta_t-s\frac{R}{\gamma-1}\Theta_x+R\frac{\Theta}{V}U_x=\kappa\left(\frac{\Theta_x}{V}\right)_x+\mu\frac{U_x^2}{V}+G,\\[4mm]
(V,U,\Theta)(0,t)=(v_-,\frac{\kappa(\gamma-1)}{\gamma
R}\frac{\Theta_{x}}{\Theta}|_{x=0}+u_b,\theta_-),\\[4mm]
(V,U,\Theta)(x,0)=(V_0,U_0,\Theta_0)=(\frac{R}{p_+}\Theta_0,\frac{\kappa(\gamma-1)}{\gamma
R}\frac{\Theta_{0x}}{\Theta_0}+u_b,\Theta_0)\to (v_+,u_b,\theta_+),\
as\ \ x\to+\infty,
\end{array}\right.\label{1.9}
\end{equation}
where
\begin{eqnarray}G&=&-\mu\frac{U_x^2}{V}=O((\ln\Theta)^2_{xx}),\nonumber\\
F(x,t)&=&\frac{\kappa(\gamma-1)}{\gamma
R}\left\{(\ln\Theta)_{xt}-s(\ln\Theta)_{xx}-\mu\left(\frac{(\ln\Theta)_{xx}}{V}\right)_x\right\}\nonumber\\
&=&\frac{\kappa a(\gamma-1)-\mu p_+\gamma}{R
\gamma}\left(\frac{(\ln\Theta)_{xx}}{\Theta}\right)_x.\label{1.10}\end{eqnarray}

Denote
\begin{eqnarray}
&&\varphi(x,t)=v(x,t)-V(x,t),\nonumber\\
&&\psi(x,t)=u(x,t)-U(x,t),\nonumber\\
&&\zeta(x,t)=\theta(x,t)-\Theta(x,t).\label{1.11}
\end{eqnarray}
Combining (\ref{1.9}) and (\ref{1.3}), the original problem can be
reformulated as
\begin{equation}\left\{
    \begin{array}{lll}
      \varphi_t-s\varphi_x=\psi_x, &  \\
      \psi_t-s\psi_x-(\frac{R\Theta}{vV}\varphi)_x+(\frac{R\zeta}{v})_x=-\mu(\frac{U_x}{vV}\varphi)_x+\mu(\frac{\psi_x}{v})_x-F, & \\
      \frac{R}{\gamma-1}\zeta_t-s\frac{R}{\gamma-1}\zeta_x+\frac{R\theta}{v}(\psi_x+U_x)-\frac{R\Theta}{V}U_x
=\kappa(\frac{\zeta_x}{v})_x-\kappa(\frac{\Theta_x\varphi}{vV})_x+\mu(\frac{{u_x}^2}{v}-\frac{{U_x}^2}{V})-G,
&\\
(\varphi,\psi,\zeta)(0,t)=(0,u_b-U(0,t),0),&\\
(\varphi,\psi,\zeta)(x,0)=(\varphi_0,\psi_0,\zeta_0)=(v_0-V_0,u_0-U_0,\theta_0-\Theta_0).&
\end{array}
  \right.
\label{1.12}\end{equation} Under the above preparation in hand, our
original problem can be transferred into a stability problem: If the
initial date $(v_0(x),u_0(x),\theta_0(x))$ of the inflow problem
(\ref{1.1}) admit a unique global solution
$(v(x,t),u(x,t),\theta(x,t))$ which tend to
$(V(x,t),U(x,t),\Theta(x,t))$ as $t\to\infty$? Recall that according
to whether $H(\R_+)$--norm of the initial perturbation
$(\varphi_0(x),\psi_0(x),\zeta_0(x))$ and (or)
$|(\theta_+-\theta_-,\rho_+-\rho_-)|>1$ or not, the stability
results are classified into global (or local) stability of strong
(or weak ) viscous contact wave.

To deduce the desired nonlinear stability result by the elementary
energy method as in
\cite{HMS2003,HMS,HMX,HXY,HZ,HLM2009,HH,TJJ2011,TZ2012}, it is
sufficient to deduce certain uniform (with respect to the time
variable $t$) energy type estimates on the solution
$(v(x,t),u(x,t),\theta(x,t))$ and the main difficulty to do so lies
in how to deal with the boundary condition when we get rid off the
small condition of $|\theta_+-\theta_-|$ and how to  establish the
Poincar$\acute{e}$ type inequality in Lemma \ref{lem3.2}  without
the smallness of $|\theta_+-\theta_-|$ which the arguments employed
in \cite{HMS2003,HMS,HMX,HXY,HZ,HLM2009,TJJ2011,TZ2012} is to use
both smallness $|\theta_+-\theta_-|$  and $N(t)=\sup_{0\leq \tau\leq
t}\|(\varphi,\psi,\zeta)\|_{H^1}$ to overcome such  difficulties.
One of the key points in such an argument is that, based on the  a
priori assumption that $\sup_{0\leq\tau\leq
t}\|(\varphi,\psi,\zeta)\|_{H^1}(\tau)$ is sufficiently small, one
can deduce a uniform lower and upper positive bounds on the specific
volume $v(x,t)$ and temperature $\theta(x,t)$. With such a bound on
$v$ and $\theta$ in hand, one can deduce certain a priori $H(\R_+)$
energy type estimates on $(\varphi,\psi,\zeta)$ in terms of the
initial perturbation $(\varphi_0,\psi_0,\zeta_0)$ provided that
$\|\Theta_{0x}\|$ suitably small, so stability of weak contact
discontinuity can be obtained . In fact if $N(t)$ not small and the
perturbation of
$\|(\varphi_{0x},\psi_{0x},\zeta_{0x})\|_{L^2(\R_+)}$ not small (see
\cite{HH}), the combination of  the  analysis similar as above with
the standard continuation argument, we can obtain the upper and
lower bounds of $(v,\theta)$, then that yields the global stability
of strong viscous contact discontinuity for the one-dimensional
compressible Navier-Stokes equations . So it is important to finish
a priori estimate without the smallness of $|\theta_+-\theta_-|$ and
$\|(\varphi_{0x},\psi_{0x},\zeta_{0x})\|_{L^2(\R_+)} .$

 This paper  we replace self-similar solution (see \cite{{HMS2003,HMS,HMX,HXY,HZ,HLM2009}}) to a diffusion equation's solution. We  use the fundamental solution skill in
\cite{TZ2012}  and give some precise time estimates about
temperature $\Theta$ which can cause to the global uniform time
estimate, so the similar energy priori estimate as the refers can be
obtained . The global stability result comes out because of these
time estimates. It is easy to see that in such a result,
  for all $t\in\R_+$, Osc
 $\theta(t):=\sup_{x\in\R_+}\theta(x,t)-\inf_{x\in\R_+}\theta(x,t)\geq|\theta_+-\theta_-|$,
 the oscillation of the  temperature $\theta(x,t)$ should not be
 sufficiently small .

To state our main result, we assume throughout of this section that
$$(\varphi_0,\zeta_0)(x)\in H_0^1(0,\infty),\ \ \psi_0(x)\in H^1(0,\infty).$$
Moreover, for an interval $I\in [0,\infty)$ , we define the function
space
$$X(I)=\left\{(\varphi,\psi,\zeta)\in C(I,H^1)|\varphi_x\in
L^2(I;L^2), (\psi_x,\zeta_x)\in L^2(I;H^1)\right\}.$$

Our  main results of this paper now reads as follows.
\begin{thm}\label{thm1.1}
There exist  positive constants $C$,$\alpha$, $\delta_0$  and
$\eta_0$ such that if $1<|\theta_+-\theta_-|$, $\delta_0$
independent of $\theta_{\pm}$ ,
$$\Theta_0=\theta_+-(\theta_+-\theta_-)\exp\{1-(1+\alpha x)^{\delta_0}\},$$
  and
   $\|(v_0-V_0,
u_0-U_0, \theta_0-\Theta_0)\|_{L^2}\leq \eta_0,$ $\|(v_{0x}-V_{0x},
u_{0x}-U_{0x}, \theta_{0x}-\Theta_{0x})\|_{L^2}\leq C,$(\ref{1.12})
has a unique global solution $(\varphi,\psi,\zeta)$ satisfying
$(\varphi,\psi,\zeta)\in X([0,\infty))$ and
 $$\sup_{x\in \R_+}|(\varphi,\psi,\zeta)|\to 0,\ as\ t\to \infty.$$
\end{thm}

{\re The constant $\alpha$ will be determined in Lemma \ref{lem2.4}
for  the definition of viscous contact discontinuity in \cite{HMS},
which is on any finite-time interval as $k \to 0$ , $(V,U,\Theta)$
is a viscous contact wave when
$\|(V-\overline{V},U-\overline{U},\Theta-\overline{\Theta})\|_{L^p}\to
0$.}

\section{Preliminary} In this section, to study the asymptotic
behavior of the solution to inflow problem (\ref{1.3}), we will do
some preparations lemmas and  list a priori estimate   which are
important to the proof of Theorem \ref{thm1.1}.

Throughout this paper, we shall denote $H^l(\R_+)$ the usual $l-th$
order Sobolev space with the norm
$$\|f\|_{l}=\big(\sum_{j=0}^l\|\partial_x^j f\|^2\big)^{1/2},\ \ \|\cdot\|:=\|\cdot\|_{L^2(\R_+)}.$$
For simplicity, we also use $C$ or $C_i$ ($i=1,2,3.....$) to denote
the various positive generic constants.  $C(z)$ stands for constant
about $z$ and $\lim_{z\to 0}C(z)=0$. $\epsilon$ and
$\epsilon_i$($i=1,2,3.....$) stand for suitably small positive
constant in Cauchy-Schwarz inequality and
$\partial_x^i=\frac{\partial^i}{\partial x^i}$.

We shall prove Theorem \ref{thm1.1} by combining the local existence
and the global-in-time a priori estimates. Since the local existence
of the solution  is well known (see, for example,\cite{HMS}), we
omit it here for brevity. to prove the global existence part of
Theorem \ref{thm1.1}, it is sufficient to establish the following a
priori estimates.

\begin{pro}\label{pro2.2}{\rm(A priori estimate)} Let $(\varphi,\psi,\zeta)\in
X([0,t])$ be a solution of problem (\ref{1.12}) for some $t>0$. Then
there exist positive  constants $C(\delta_0)<1$ and $C$ which are
all independent of $t$ and $(v,\theta)$ , such that if $m\leq
v,\theta\leq M$
 and $N(t)=\sup_{0\leq \tau\leq t}\|(\varphi,\psi,\zeta)\|_1\leq
C$,
 it holds that
\begin{eqnarray}\label{2.1}
&&\sup_{0\leq \tau\leq t}\|(\psi,\varphi,\zeta)\|^2(t)+\int_0^t\|(\psi_x,\zeta_x)\|^2(\tau)d\tau\nonumber\\
&&\leq C\|(\varphi_0,\psi_0,\zeta_0)\|+C(\delta_0)
.\nonumber\\
&&\sup_{0\leq \tau\leq t}\|(\psi_x,\varphi_x,\zeta_x)\|^2(t)+\int_0^t\left(\|\varphi_x\|^2(\tau)+\|(\psi_x,\zeta_x)\|^2_1(\tau)\right)d\tau\nonumber\\
&&\leq C\|(\varphi_0,\psi_0,\zeta_0)\|_1+C(\delta_0).\end{eqnarray}
\end{pro}

To finish this proposition,  we must consider  some properties of
$\Theta_0$ and $\partial_x^i\Theta$ $(i=1,2,3...)$ as we list
following.
\begin{lem}\label{lem2.1}
As to the definition of $\Theta_0$ in Theorem \ref{thm1.1} we have
\begin{eqnarray*}
&&\|\Theta_0-\theta_+\|_{L^1}\leq C\alpha^{-1}\sum_{n=0}^{[\frac{1}{\delta_0}]-1}\prod_{i=0}^n(\frac{1}{\delta_0}-i),\\
&&0<\Theta_{0x}\leq C\alpha\delta_0(1+\alpha x)^{\delta_0-1}\exp\{-(1+\alpha x)^{\delta_0}\},\\
&&|\Theta_{0xx}|\leq C\alpha^2\delta_0\big((1+\alpha x)^{2\delta_0-2}+(1-\delta_0)(1+\alpha x)^{\delta_0-2}\big)\exp\{-(1+\alpha x)^{\delta_0}\},\\
&&\|\Theta_{0x}\|^2\leq C\alpha\delta_0,\\
&&\|\Theta_{0x}\|_{L^1(\R_+)}\leq C,\\
&&\|\Theta_{0xx}\|^2+\|(\ln\Theta_0)_{xx}\|^2\leq
C\alpha^3\delta_0^2,\\
&&\|\Theta_{0xxx}\|^2+\|(\ln\Theta_0)_{xxx}\|^2\leq
C,\\
&&\int_{\R}\Theta_{0x}^2(1+\alpha x)dx\leq C\alpha\delta_0.
\end{eqnarray*}
\end{lem}
\pf In fact
\begin{eqnarray*}
\int_{\R_+}|\Theta_0-\theta_+|dx&\leq&
C\alpha^{-1}\int_{\R_+}\exp\{-(1+\alpha x)^{\delta_0}\}d\left((\alpha x+1)^{\delta_0}\right)^{1/\delta_0}\\
&=&C\alpha^{-1}\int_0^{+\infty}\exp\{-(1+\alpha x)^{\delta_0}\}\frac{1}{\delta_0}\left((1+\alpha x)^{\delta_0}\right)^{1/\delta_0-1}d(\alpha x+1)^{\delta_0}\\
&=&C\alpha^{-1}\frac{1}{\delta_0}\int_1^{+\infty}\exp\{-z\}z^{1/\delta_0-1}dz\\
&=&C\alpha^{-1}\frac{1}{\delta_0}\exp\{-z\}z^{1/\delta_0-1}\Big|^1_{+\infty}+C\alpha^{-1}\frac{1}{\delta_0}\int_1^{+\infty}\exp\{-z\}dz^{1/\delta_0-1}\\
&=&C\alpha^{-1}\frac{1}{\delta_0}\exp\{-z\}z^{1/\delta_0-1}\Big|^1_{+\infty}+...
+C\alpha^{-1}\prod_{i=0}^{[\frac{1}{\delta_0}]-1}(\frac{1}{\delta_0}-i)\int_1^{+\infty}\exp\{-z\}z^{1/\delta_0-[1/\delta_0]}dz\\
&\leq&
C\alpha^{-1}\sum_{n=0}^{[\frac{1}{\delta_0}]-1}\prod_{i=0}^n(\frac{1}{\delta_0}-i).
\end{eqnarray*}
That is the first inequality.

 As to the last inequality, we know that
\begin{eqnarray*}
\int_{\R}|\Theta_{0x}|^2(1+\alpha x)dx&\leq&\int_{\R}C\alpha^2\delta_0^2(1+\alpha x)^{2\delta_0-2}(1+\alpha x)\exp\{-2(1+\alpha x)^{\delta_0}\}dx\\
&=&\int_{\R}C\alpha^2\delta_0^2(1+\alpha x)^{2\delta_0-1}\exp\{-2(1+\alpha x)^{\delta_0}\}dx\\
&=&\int_{\R}C\alpha\delta_0^2(2\delta_0)^{-1}\exp\{-2(1+\alpha x)^{\delta_0}\}d(1+\alpha x)^{2\delta_0}\\
&=&\int_{\R}C2^{-1}\alpha\delta_0\exp\{-2z\}dz^2\leq
C\alpha\delta_0.
\end{eqnarray*}

The other  inequalities can be check easily by using the definition
of $\Theta_0$ and we omit to write them. $\Box$

 Next,  we construct a parabolic equation about
 $\theta_2$, it will be used in the estimates of  $\partial_x^i\Theta$.

\begin{lem}\label{lem2.2}If $\delta_0$ and $\Theta_0$ satisfying the
condition in Theorem \ref{thm1.1} and
\begin{eqnarray*}&&\theta_2(x,t)=\int_{0}^{+\infty}(4\pi
at)^{-1/2}(\Theta_0(h)-\theta_-)\left\{\exp\{-\frac{(h-st-x)^2}{4at}\}-\exp\{-\frac{(h-st+x)^2}{4at}\}\right\}\
dh+\theta_-,\\
&&K=\int_0^{+\infty}(4\pi
at)^{-1/2}\Theta_{0z}(z)\exp\{-\frac{(z-st+x)^2}{4at}\}dz,\end{eqnarray*}
 we can get
 \begin{eqnarray}
 &&\theta_{2t}-s\theta_{2x}=a\theta_{2xx}-2sK;\nonumber\\
 &&\theta_2(0,t)=\theta_-;\nonumber\\
 &&\theta_2(x,0)=\theta_{20}(x)=\left\{
                      \begin{array}{ll}
                        \Theta_{0}(x)\to \theta_+ , & \hbox{$x>0$;} \\
                       -\Theta_{0}(-x)+2\theta_-\to 2\theta_--\theta_+ , & \hbox{$ x\leq 0$,}
                      \end{array}
                    \right.\label{2.2}
 \end{eqnarray}
 and\begin{eqnarray}&&\int_0^t\|K\|_{L^1(\R_+)}dt\leq C\alpha^{1/2}\delta_0^{1/2}(1+t)^{1/2},\label{2.3}\\
 &&\int_0^t\|\theta_{2x}\|^2dt\leq C(1+t)^{1/2},\label{2.4}\\
  &&\|\theta_{2x}\|^2+\int_0^t\ \theta_{2x}^2(0,t)\ dt+\int_0^t\ \int_0^{+\infty}\ \theta_{2xx}^2\ dx\ dt\leq
 C\alpha\delta_0.\label{2.5}\end{eqnarray}

\end{lem}
\pf First we proof (\ref{2.2})$_1$ as following:

Because
\begin{eqnarray}\label{2.7}
\theta_{2x} &=&\int_0^{+\infty}(4\pi
at)^{-1/2}\left(\Theta_0(z)-\theta_-\right)\exp\{-\frac{(z-st-x)^2}{4at}\}\frac{z-st-x}{2at}dz\nonumber\\
&&\quad+\int_0^{+\infty}(4\pi
at)^{-1/2}\left(\Theta_0(z)-\theta_-\right)\exp\{-\frac{(z-st+x)^2}{4at}\}\frac{z-st+x}{2at}dz\nonumber\\
 &=&\int_0^{+\infty}(4\pi
at)^{-1/2}\Theta_{0z}(z)\exp\{\frac{-(z-st-x)^2}{4at}\}dz\nonumber\\
&&\quad-\int_0^{+\infty}(4\pi
at)^{-1/2}\Theta_{0z}(z)\exp\{\frac{-(z-st+x)^2}{4at}\}dz\nonumber\\
&=&:\hat{I}_1+\hat{I}_2 ,\end{eqnarray}  it is easy to check
$\theta_{2t}-s\theta_{2x}=a\theta_{2xx}-2sK$, $K=-\hat{I}_2$ and we
finish (\ref{2.2})$_1$ and (\ref{2.2})$_2$.

Now we proof (\ref{2.2})$_3$ as following.

 From heat conduction equation's  initial theorem  and uniform estimates we know that
\begin{eqnarray*}
\theta_{20}(x)&=&\lim_{t\to 0}\int_{-\infty}^{+\infty}(4\pi
at)^{-1/2}\theta_{20}(h)\exp\{-\frac{(h-x)^2}{4at}\}dh\\
 &=&\lim_{t\to
0}\int_0^{+\infty}(4\pi
at)^{-1/2}(\Theta_0(h)-\theta_-)\exp\{-\frac{(h-x)^2}{4at}\}dh\\
&&+\lim_{t\to 0}\int_0^{+\infty}(4\pi
at)^{-1/2}(-\Theta_0(h)+\theta_-)\exp\{-\frac{(h+x)^2}{4at}\}dh+\theta_-.
\end{eqnarray*}
So
\begin{eqnarray*}
\lim_{t\to 0}(\theta_2(x,t)-\theta_{20}(x))&=&\lim_{t\to
0}\int_0^{+\infty}(4\pi
at)^{-1/2}(\Theta_0(h)-\theta_-)\left(\exp\{-\frac{(h-x-st)^2}{4at}\}-\exp\{-\frac{(h-x)^2}{4at}\}\right)dh\\
&&+\lim_{t\to 0}\int_0^{+\infty}(4\pi
at)^{-1/2}(-\Theta_0(h)+\theta_-)\left(\exp\{\frac{-(h+x-st)^2}{4at}\}-\exp\{-\frac{(h+x)^2}{4at}\}\right)dh.
\end{eqnarray*}
Use Lebesgue control theorem we know
\begin{eqnarray*}
&&\lim_{t\to 0}\int_0^{+\infty}(4\pi
at)^{-1/2}|\Theta_0(h)-\theta_-|\left|\exp\{-\frac{(h-x-st)^2}{4at}\}-\exp\{-\frac{(h-x)^2}{4at}\}\right|dh\\
&&+\lim_{t\to 0}\int_0^{+\infty}(4\pi
at)^{-1/2}|-\Theta_0(h)+\theta_-|\left|\exp\{\frac{-(h+x-st)^2}{4at}\}-\exp\{-\frac{(h+x)^2}{4at}\}\right|dh\\
&&\leq\lim_{t\to
0}C\int_0^{+\infty}\left|\exp\{-\frac{(h-x-st)^2}{4at}\}-\exp\{-\frac{(h-x)^2}{4at}\}\right|d
(4
at)^{-1/2}h\\
&&+\lim_{t\to
0}C\int_0^{+\infty}\left|\exp\{\frac{-(h+x-st)^2}{4at}\}-\exp\{-\frac{(h+x)^2}{4at}\}\right|d
(4
at)^{-1/2}h\\
&&\leq C\int_{-\infty}^{+\infty}e^{-\xi^2}\lim_{t\to
0}\left(\exp\{-(\xi-\frac{st}{\sqrt{4at}})^2+\xi^2\}-1\right)d\xi=0,
\end{eqnarray*}
which means $$\lim_{t\to 0}(\theta_2(x,t)-\theta_{20}(x))=0.$$ So
 $(\ref{2.2})_3$ is established.

Now we consider the estimate about $K$.  In fact from (\ref{2.7}) we
know
\begin{eqnarray}\label{2.9}
K=-\hat{I}_2&=&\int_0^{+\infty}(4\pi
at)^{-1/2}\Theta_{0z}(z)\exp\{-\frac{(z-st+y)^2}{4at}\}dz,
\end{eqnarray}
so we can get
\begin{eqnarray*}
&&\int_0^t\int_0^{+\infty}|K|dxdt=
\int_0^t\int_0^{+\infty}|\hat{I}_2|dxdt\\
&&\leq C\int_0^t(4\pi
at)^{-1/2}\|\Theta_{0z}\|\left(\int_0^{+\infty}\exp\{-\frac{(z+x)^2}{4at}\}\exp\{-\frac{s^2t}{4a}\}dz\right)^{1/2}dt\\
&&\quad\times\int_0^{+\infty}(4\pi
at)^{-1/2}\exp\{-\frac{(z+x)^2}{4at}\}dx\\
&&\leq C\|\Theta_{0z}\|(1+t)^{1/2}\leq
C\alpha^{1/2}\delta_0^{1/2}(1+t)^{1/2}.\end{eqnarray*}

And use H$\ddot{o}$lder inequality and Fubini Theorem and
$1<|\theta_+-\theta_-|\leq\|\Theta_{0z}\|_{L^1(\R_+)}<C$ (see Lemma
\ref{lem2.1}), we can get
\begin{eqnarray*}\int_0^t\int_0^{+\infty}|\hat{I}_1|^2dxdt&\leq&
C\int_0^t\int_0^{+\infty}(4\pi
at)^{-1}\left(\int_0^{+\infty}\Theta_{0z}\exp\{-\frac{(z-st-x)^2}{4at}\} dz\right)^2dxdt\\
&\leq&C\int_0^t\int_0^{+\infty}(4\pi
at)^{-1}\int_0^{+\infty}|\Theta_{0z}|\exp\{-\frac{(z-st-x)^2}{4at}\}
dzdx\int_0^{+\infty}|\Theta_{0z}|dzdt\\
&\leq &C\sqrt{1+t}.\end{eqnarray*}

In all when we combine with the estimates about $\hat{I}_{1}$ and
$\hat{I}_{2}$ of (\ref{2.7}) we can get
$$\int_0^t\int_0^\infty \theta_{2x}^2dxdt\leq
C\sqrt{1+t}.$$

Now both side of (\ref{2.2})$_1$ multiply by $\theta_{2xx}$,
integrate in $\R_+\times(0,t)$ and combine with Cauchy-Schwarz
inequality we can get \begin{equation}\|\theta_{2x}\|^2+\int_0^t\
\theta_{2x}^2(0,t)\ dt+\int_0^t\ \int_0^{+\infty}\ \theta_{2xx}^2\
dx\ dt\leq C\|\Theta_{0x}\|^2\leq C\alpha
\delta_0.\label{2.10}\end{equation} So we finish this lemma.$\Box$

 Now let's consider the time estimates about $\partial_x^i\Theta$ $(i=1,2,3)$ of
 (\ref{1.8}), we have the following results. We list the proof steps of each
formula
 inside this lemma for reading convenient.
\begin{lem}\label{lem2.3}
If $\Theta_{0x}$ satisfying the condition of Theorem \ref{thm1.1},we
can get
\begin{eqnarray}&&\|(\ln\Theta)_x\|^2+\int_0^t\
(\ln\Theta)^2_x(0,t)\ dx+a\int_0^t\ \|(\ln\Theta)_{xx}\|^2\ dt \leq
C\alpha\delta_0.\label{2.11}\\
&&_{(see (\ref{2.21})-(\ref{2.23}))}\nonumber\\
 &&\|\Theta-\theta_2\|^2+\int_0^t\ \|(\ln\Theta)_{x}\|^2\ dt\leq
C(1+t)^{1/2}.\label{2.12}\\
&&_{{(see (\ref{2.19})-(\ref{2.24}))}}\nonumber\\
 && \|(\ln\Theta)_x\|^2\leq
C(1+t)^{-1/2}.\label{2.13}\\&&_{{(see
(\ref{2.25})-(\ref{2.28}))}}\nonumber\\
&&\|(\ln\Theta)_{xx}\|^2\leq C(1+t)^{-3/2}.\label{2.14}\\.&&_{(see
(\ref{2.29})-(\ref{2.35}))}\nonumber\\
&&\|(\ln\Theta)_{xx}\|^2(1+t)+\int_0^t\|\partial^3_x\ln\Theta\|^2(1+t)\
dt\nonumber\\ &&\quad+\int_0^t\ (\partial^2_x\ln\Theta)^2(0,t)(1+t)\
dt\leq
C\delta_0.\label{2.15}\\&&_{(see(\ref{2.36}))}\nonumber\\
&& \|\partial^3_x\ln\Theta\|^2\leq
C(1+t)^{-5/2}.\label{2.16}\\&&_{(see
(\ref{2.37})-(\ref{2.41}))}\nonumber\\
 &&\int_0^t(\partial^3_x\ln\Theta)^2(0,t)\ dt\leq C.\label{2.17}\\&&_{(see
(\ref{2.42}))} \nonumber\\
&&\int_{\R_+}\Theta_x^2xdx\leq C\delta_0.\label{2.18}\\
&&_{(see (\ref{2.43})-(\ref{2.44}))} \nonumber
\end{eqnarray}
\end{lem}
\pf

Both side of (\ref{1.8})$_1-$(\ref{2.2})$_1$ multiply by
$\Theta-\theta_2$, integrate in $\R_+\times (0,t)$ and combine with
Cauchy-Schwarz inequality we can get
\begin{eqnarray}&&\|\Theta-\theta_2\|^2+\int_0^t\|(\ln\Theta)_x\|^2\
dt\nonumber\\
&& \leq C\int_0^t(\|{\theta_2}_x\|^2+\|K\|_{L^1})\ dt+C\int_0^t\
\left((\ln\Theta)^2_{x}+{\theta}^2_{2x}\right)(0,t)\
dt.\label{2.19}\end{eqnarray}

On the other hand from (\ref{1.8})$_1$
\begin{equation}(\ln\Theta)_t-s(\ln\Theta)_x=a\frac{(\ln\Theta)_{xx}}{\Theta},\label{2.20}\end{equation}
both side of it multiply by $(\ln\Theta)_{xx}$ and integrate in
$\R_+\times (0,t)$ we can get

\begin{eqnarray}&&\|(\ln\Theta)_x\|^2+\int_0^t\
(\ln\Theta)^2_x(0,t)\ dx+a\int_0^t\ \|(\ln\Theta)_{xx}\|^2\
dt\nonumber\\
&&\leq
C\|(\ln\Theta_{0})_x\|^2+\int_0^t(\ln\Theta)_t(\ln\Theta)_x\big|_0^{+\infty}\
dt.\label{2.21}\end{eqnarray}

According to Lemma \ref{lem2.1} we know that
$\|(\ln\Theta_{0})_x\|^2\leq C\alpha\delta_0,$ when combine with
(\ref{2.21}) and
\begin{equation}\label{2.22}\Theta_t(+\infty,t)=0,\ \ \Theta_t(0,t)=0\end{equation} we can get
\begin{equation}\|(\ln\Theta)_x\|^2+\int_0^t\
(\ln\Theta)^2_x(0,t)\ dx+a\int_0^t\ \|(\ln\Theta)_{xx}\|^2\ dt\leq
C\alpha\delta_0.\label{2.23}\end{equation}

Use (\ref{2.10}), (\ref{2.3}),(\ref{2.4}) and (\ref{2.23}) to
(\ref{2.19}) we can get
\begin{equation}\|\Theta-\theta_2\|^2+\int_0^t\ \|(\ln\Theta)_{x}\|^2\ dt\leq
C(1+t)^{1/2}. \label{2.24}\end{equation} That is (\ref{2.12}).

 Next, both side of(\ref{1.8})$_1$ multiply
$\Theta^{-1}(\ln\Theta)_{xx}(1+t)$ and integrate in
$\R_+\times(0,t)$, we can get
\begin{eqnarray}&&\int_0^t(1+t)(\ln\Theta)_t(\ln\Theta)_x(0,t)\ dt+s/2\int_0^t(\ln\Theta)^2_x(0,t)(1+t)\ dt\nonumber\\
&&=a\int_0^t\int_0^{+\infty}\
\frac{(\ln\Theta)^2_{xx}}{\Theta}(1+t)\
dxdt+\int_0^t\int_0^{+\infty}\ \left((\ln\Theta)^2_x\right)_t(1+t)\
dxdt.\label{2.25}\end{eqnarray} Because
\begin{equation}\label{2.26}\int_0^t(1+t)(\ln\Theta)_t(\ln\Theta)_x(0,t)\
dt=0,\end{equation} we can get
\begin{eqnarray}
&&(1+t)\|(\ln\Theta)_x\|^2+\int_0^t\
\int_0^{+\infty}(1+t)(\ln\Theta)^2_{xx}\ dx\ dt+\int_0^t
(\ln\Theta)_x^2(0,t)(1+t)
dt\nonumber\\
 &&\leq C\|\Theta_{0x}\|^2+\int_0^t\ \int_0^{+\infty}\ (\ln\Theta)^2_x\ dx\ dt.\label{2.27}\end{eqnarray}
Combine with (\ref{2.24}) we can get
\begin{eqnarray}
&&(1+t)\|(\ln\Theta)_x\|^2+\int_0^t\
\int_0^{+\infty}(1+t)(\ln\Theta)^2_{xx}\ dx\ dt+\int_0^t
(\ln\Theta)_x^2(0,t)(1+t)
dt\nonumber\\
 &&\leq C(1+t)^{1/2}.\label{2.28}\end{eqnarray}
That means $\|(\ln\Theta)_x\|^2\leq C(1+t)^{-1/2}$, which
is(\ref{2.13}).

Again from (\ref{1.8})$_1$ we can get
\begin{equation}(\ln\Theta)_{xt}-s(\ln\Theta)_{xx}=a\left(\frac{(\ln\Theta)_{xx}}{\Theta}\right)_x.\label{2.29}
\end{equation}
Both side of (\ref{2.29})multiply $\partial^3_x\ln\Theta$ and get
\begin{equation}\left((\ln\Theta)_{xt}(\ln\Theta)_{xx}\right)_x-1/2(\partial^2_x\ln\Theta)_t
-s/2(\partial^2_x\ln\Theta)_x=a\left(\frac{(\ln\Theta)_{xx}}{\Theta}\right)_x\partial^3_x(\ln\Theta).\label{2.30}\end{equation}
Take  $\Theta_t-s\Theta_x=(a\ln\Theta)_{xx}$ into
$\left((\ln\Theta)_{xt}(\ln\Theta)_{xx}\right)_x$ we can get
\begin{eqnarray*}&&\left((\ln\Theta)_{xt}(\ln\Theta)_{xx}\right)_x(1+t)^2\\
&&=a^{-1}\left((\ln\Theta)_{xt}(\Theta_t-s\Theta_x)\right)_x(1+t)^2\\
&&=a^{-1}\left((\ln\Theta)_{xt}(\Theta_t)\right)_x(1+t)^2-sa^{-1}\left((\ln\Theta)_{xt}\Theta_x\right)_x(1+t)^2\\
&&=a^{-1}\left((\ln\Theta)_{xt}(\Theta_t)\right)_x(1+t)^2-sa^{-1}\left((\frac{\Theta_{xt}}{\Theta}-\frac{\Theta_t\Theta_x}{\Theta^2})\Theta_x\right)_x(1+t)^2\\
&&=a^{-1}\left((\ln\Theta)_{xt}(\Theta_t)\right)_x(1+t)^2-sa^{-1}\left(\frac{1}{2}(\frac{\Theta^2_x}{\Theta})_t
+\frac{1}{2}\frac{\Theta^2_x\Theta_t}{\Theta^2}-\frac{\Theta_t\Theta^2_x}{\Theta^2}\right)_x(1+t)^2\\
&&=a^{-1}\left((\ln\Theta)_{xt}(\Theta_t)\right)_x(1+t)^2-sa^{-1}/2\left((\frac{\Theta^2_x}{\Theta}(1+t)^2)_{tx}-2(\frac{\Theta^2_x}{\Theta})_x(1+t)-(\frac{\Theta^2_x\Theta_t}{\Theta^2})_x(1+t)^2\right).
\end{eqnarray*}
Now both side of (\ref{2.30}) multiply $(1+t)^2$ then integrate in
$\R_+\times(0,t)$ and combine with (\ref{2.22}),
$\Theta_x(\infty,t)=0$ and Cauchy-Schwarz inequality we get for a
small $\epsilon>0$ ,
\begin{eqnarray}
&&-s/(2a)\int_0^{+\infty}\
\left(\frac{\Theta_x^2}{\Theta}(1+t)^2\right)_xdx+s/(2a)\int_0^{+\infty}\
\left(\frac{\Theta_{0x}^2}{\Theta_0}\right)_x\
dx-s/a\int_0^t\frac{\Theta_x^2}{\Theta}(0,t)(1+t)\ dt\nonumber\\
&&\geq-s/2\int_0^t\ (\ln\Theta)_{xx}^2(0,t)(1+t)^2\ dt+a\int_0^t\
\int_0^{+\infty}\ \frac{(\ln\Theta)_{xxx}^2}{\Theta}(1+t)^2\ dx\
dt\nonumber\\
&&\quad-\epsilon\int_0^t\ \int_0^{+\infty}\
(1+t)^2(\ln\Theta)_{xxx}^2\ dx\ dt-Ca\epsilon^{-1}\int_0^t\
\int_0^{+\infty}\ (1+t)^2(\ln\Theta)_{xx}^2(\ln\Theta)_x^2\ dx\
dt\nonumber\\
&&\quad+1/2\|(\ln\Theta)_{xx}\|^2(1+t)^2-1/2\|(\ln\Theta_0)_{xx}\|^2-\int_0^t\
\|(\ln\Theta)_{xx}\|^2(1+t)\ dx\nonumber\\
&&\geq-s/2\int_0^t\ (\ln\Theta)_{xx}^2(0,t)(1+t)^2\ dt+Ca\int_0^t\
\int_0^{+\infty}\ \frac{(\ln\Theta)_{xxx}^2}{\Theta}(1+t)^2\ dx\
dt\nonumber\\
&&\quad-C\epsilon^{-1}a\int_0^t\ \int_0^{+\infty}\
(1+t)^2\|(\ln\Theta)_{xx}\|\|(\ln\Theta)_{xxx}\|(\ln\Theta)_x^2\ dx\
dt\nonumber\\
&&\quad+1/2\|(\ln\Theta)_{xx}\|^2(1+t)^2-1/2\|(\ln\Theta_0)_{xx}\|^2-\int_0^t\
\|(\ln\Theta)_{xx}\|^2(1+t)\ dx.\label{2.31}
\end{eqnarray}
From Lemma \ref{lem2.1}we know that
\begin{equation}\frac{\Theta_{0x}^2}{\Theta_0}(0)\leq C\alpha^2\delta_0^2,\label{2.32}\end{equation}
$$\|(\ln\Theta_0)_{xx}\|^2\leq C\alpha^3\delta_0^2.$$

Combine with (\ref{2.28}) we can get
\begin{equation}\left|s/(2a)\int_0^{+\infty}\
\left(\frac{\Theta_{0x}^2}{\Theta_0}\right)_x\
dx-s/a\int_0^t\frac{\Theta_x^2}{\Theta}(0,t)(1+t)\ dt\right|\leq
C(1+t)^{1/2}.\label{2.33}\end{equation} Take (\ref{2.28})
(\ref{2.32}) and (\ref{2.33}) into (\ref{2.31}) we can get

\begin{eqnarray}&&\|(\ln\Theta)_{xx}\|^2(1+t)^2+\int_0^{t}\
(1+t)^2(\ln\Theta)_{xx}^2(0,t)\ dt+\int_0^t\ \int_0^{+\infty}\
(1+t)^2(\ln\Theta)_{xxx}^2\ dx\ dt\nonumber\\
&& \leq C(1+t)^{1/2},\label{2.34}\end{eqnarray} which also means

\begin{equation}\|(\ln\Theta)_{xx}\|^2\leq
C(1+t)^{-3/2},\label{2.35}\end{equation} and finish (\ref{2.14}).

If both side of (\ref{2.30}) multiply by $(1+t)$, similar as the
proof of (\ref{2.34}),  when combine with (\ref{2.23}) we can get
\begin{equation}\|(\ln\Theta)_{xx}\|^2(1+t)+\int_0^t\ \int_0^{+\infty}\ (1+t)(\partial^3_x\ln\Theta)^2\ dx\
dt+\int_0^t\ (\partial^2_x\ln\Theta)^2(0,t)(1+t)\ dt\leq
C\delta_0,\label{2.36}\end{equation} which means (\ref{2.15}).

From (\ref{2.29}) we can get
\begin{equation}
\partial_t(\ln\Theta)_{xx}-s\partial_x(\ln\Theta)_{xx}=a\partial^2_x\left(\frac{(\ln\Theta)_{xx}}{\Theta}\right).\label{2.37}
\end{equation} Similar as (\ref{2.30}) we need to deal with the
boundary term about
$\left((\ln\Theta)_{xxt}(\ln\Theta)_{xxx}\right)_x$.

Because combine with (\ref{2.20}) and (\ref{2.29})we can get
\begin{eqnarray*}
&&\left((\ln\Theta)_{xxt}(\ln\Theta)_{xxx}\right)_x\\
&&=1/s\left((\ln\Theta)_{xxt}(\ln\Theta)_{xx}(\ln\Theta)_t\right)_x-a/s\left((\ln\Theta)_{xxt}\frac{(\ln\Theta)_{xx}^2}{\Theta}\right)_x\\
&&\quad+1/a\left((\ln\Theta)_{xxt}\Theta(\ln\Theta)_{xt}\right)_x-s/a\left((\ln\Theta)_{xxt}\Theta(\ln\Theta)_{xx}\right)_x\\
&&:=I_1+I_2+I_3+I_4,
\end{eqnarray*}
when both side of (\ref{2.37}) multiply
$\partial^4_x\ln\Theta(1+t)^3$ then integrate in $\R_+\times(0,t)$
we can get
\begin{eqnarray}
&&\int_0^t\int_0^{+\infty}\left(\big((\ln\Theta)_{xxt}(\ln\Theta)_{xxx}\big)_x-s/2\partial_x(\ln\Theta)_{xx}^2\right)(1+t)^3dxdt\nonumber\\
&&=\int_0^t\int_0^{+\infty}\left(I_1+I_2+I_3+I_4-s/2\partial_x(\ln\Theta)_{xxx}^2\right)(1+t)^3dxdt\nonumber\\
&&=\int_0^t\int_0^{+\infty}a\partial^2_x\left(\frac{(\ln\Theta)_{xx}}{\Theta}\right)\partial^4_x\ln\Theta(1+t)^3dxdt\nonumber\\
&&\quad+\int_0^t\int_0^{+\infty}\frac{1}{2}\left((\partial^3_x\ln\Theta)^2\right)_t(1+t)^3dxdt.\label{2.38}
\end{eqnarray}

To finish (\ref{2.38}),similar as (\ref{2.34}) using (\ref{2.28})
and (\ref{2.34}) we can get
\begin{equation}\|\partial^3_x\ln\Theta\|^2(1+t)^3+\int_0^t\
(1+t)^3\|\partial^4_x\ln\Theta\|^2\ dt+\int_0^t\
(\partial^3_x\ln\Theta)^2(1+t)^3(0,t)\ dt\leq
C(1+t)^{1/2}.\label{2.41}\end{equation} This means (\ref{2.16})
finished.

 When we change $(1+t)^3$  to
$(1+t)^2$ and combine with (\ref{2.15}), we can get

\begin{equation}\|\partial^3_x\ln\Theta\|^2(1+t)^2+\int_0^t\
(1+t)^2\|\partial^4_x\ln\Theta\|^2\ dt+\int_0^t\
(\partial^3_x\ln\Theta)^2(1+t)^2(0,t)\ dt\leq
C,\label{2.42}\end{equation} which finish (\ref{2.17}).

Now both side of (\ref{2.20}) multiply by
$(\ln\Theta)_{xx}(2s(\tau-t)+x)$ and integrate in
$[-2s(\tau-t),\infty)\times (0,t)$ , we can get
\begin{eqnarray}\label{2.43}
\sum_{i=1}^6K_i&=&\int_0^t\int_{-2s(\tau-t)}^{+\infty}\left((\ln\Theta)_\tau(\ln\Theta)_x(2s(\tau-t)+x)\right)_xdxd\tau-1/2\int_0^t\int_{-2s(\tau-t)}^{+\infty}\left((\ln\Theta)_x^2(2s(\tau-t)+x)\right)_\tau dxd\tau\nonumber\\
&&\quad-a\int_0^t\int_{-2s(\tau-t)}^{+\infty}(\ln\Theta)_{xx}(\ln\Theta)_x\Theta^{-1}dxd\tau-s/2\int_0^t\int_{-2s(\tau-t)}^{+\infty}\left((\ln\Theta)_x^2(x+2s(\tau-t))\right)_xdxd\tau\nonumber\\
&&\quad+\frac{s}{2}\int_0^t\int_{-2s(\tau-t)}^{+\infty}(\ln\Theta)_x^2dxd\tau-\int_0^t\int_{-2s(\tau-t)}^{+\infty}a(\ln\Theta)_{xx}^2(x+2s(\tau-t))\Theta^{-1}
dxd\tau=0.
\end{eqnarray}
Use Cauchy-Schwarz inequality
\begin{eqnarray*}
|K_3|\leq
\frac{-s}{4}\int_0^t\int_{-2s(\tau-t)}^{+\infty}(\ln\Theta)_x^2dxd\tau+C\int_0^t\int_{-2s(\tau-t)}^{+\infty}(\ln\Theta)_{xx}^2dxd\tau,
\end{eqnarray*}
then combine with (\ref{2.11}) and $s<0$ we can get
$$|K_3|+K_5\leq C\int_0^t\int_{-2s(\tau-t)}^{+\infty}(\ln\Theta)_{xx}^2dxd\tau+\frac{s}{4}\int_0^t\int_{-2s(\tau-t)}^{+\infty}(\ln\Theta)_x^2dxd\tau\leq C\delta_0+\frac{s}{4}\int_0^t\int_{-2s(\tau-t)}^{+\infty}(\ln\Theta)_x^2dxd\tau.$$

Use parabolic extreme value theory we know $|(\ln\Theta)_xx|\leq C$,
so $\lim_{x\to+\infty}\Theta_x^2(x+2s(\tau-t))=0$. When combine with
the estimates from $K_1$ to $K_5$, Lemma \ref{lem2.1} and
$\lim_{x\to +\infty}\Theta_x^2(x+2s(\tau-t))=0$, (\ref{2.43}) can be
change to
\begin{eqnarray}\label{2.44}
&&\int_0^{+\infty}(\ln\Theta)_x^2xdx+\int_0^t\int_{-2s(\tau-t)}^{+\infty}(\ln\Theta)_{xx}^2(x+2s(\tau-t))dxd\tau\nonumber\\
&&\quad-\frac{s}{4}\int_0^t\int_{-2s(\tau-t)}^{+\infty}(\ln\Theta)_x^2dxd\tau
\leq C\delta_0.
\end{eqnarray}
 So we finish this lemma. $\Box$

The next lemma is concerned with the relations between the viscous
continuity and the contact discontinuity. We shall show that as the
heat conductivity $k$ goes to zero, $(V,U,\Theta)$ will approximate
$(\overline{V},\overline{U},\overline{\Theta})$ in $L^p(\R_+)$
$(p\geq 1)$ norm on any finite time interval.

\begin{lem}\label{lem2.4}
For any given $T\in(0,+\infty)$ independent of $\kappa$ such that
for any $p\geq 1$ and $t\in[0,T]$,
$$\|(V-\overline{V},U-\overline{U},\Theta-\overline{\Theta})\|_{L^p(\R_+)}\to 0,\ \ \mathrm{as}\ \ \kappa\to 0.$$
\end{lem}
\pf

Letting $$\Omega_1=(0,-st)\ \ \mathrm{and} \ \
\Omega_2=(-st,+\infty).$$ By the definition of $\overline{\Theta}$
in (\ref{1.5}), to estimate
$\|\Theta-\overline{\Theta}\|_{L^p(\R_+)}$, it suffices to prove

$$\|\Theta-\theta_-\|_{L^p(\Omega_1)},\ \ \|\Theta-\theta_+\|_{L^p(\Omega_2)}\to 0,\ \ \mathrm{as}\ \ \kappa\to 0,\ p\geq1.$$
Because$$\|\Theta-\theta_-\|^p_{L^p(\Omega_1)}\leq
C\|\Theta-\theta_-\|_{L^1(\Omega_1)},\ \
\|\Theta-\theta_+\|^p_{L^p(\Omega_2)}\leq
C\|\Theta-\theta_+\|_{L^1(\Omega_2)},$$ the only thing we need to
proof is
$$\lim_{\kappa\to 0}\|\Theta-\theta_-\|_{L^1(\Omega_1)}+\|\Theta-\theta_+\|_{L^1(\Omega_2)}=0.$$

 In fact we set $sgn_{\eta}(l)=\left\{
                   \begin{array}{ll}
                     1, & \hbox{$l> \eta$;} \\
                     l/\eta, & \hbox{$-\eta\leq l\leq\eta$;} \\
                     -1, & \hbox{$l< -\eta$.}
                   \end{array}
                 \right.
$, $I_{\eta}(l)=\int_0^l sgn_{\eta}(l) dl$ and $\eta>0$. Both side
of (\ref{1.8})$_1$ multiply by $sgn_{\eta}(\Theta-\theta_-)$ and
integrate in $(0,-s\tau)\times (0,t)$ we can get
\begin{eqnarray*}
\int_0^t\left(\int_0^{-s\tau}I_{\eta}(\Theta-\theta_-)dx\right)_{\tau}
d\tau&=&a\int_0^t(\ln\Theta)_x(-st,t)sgn_{\eta}(\Theta-\theta_-)(-st,t)d\tau\\
&&\quad-a\int_0^t\int_0^{-s\tau}(\ln\Theta)_x^2sgn'_{\eta}(\Theta-\theta_-)dxd\tau.
\end{eqnarray*}
So  when $\eta\to 0$ and use (\ref{2.11}) and (\ref{2.15}) we can
get
\begin{eqnarray}\label{2.46}
&&\|\Theta-\theta_-\|_{L^1(\Omega_1)}
+a\int_0^t\int_0^{-s\tau}(\ln\Theta)_x^2sgn'_{\eta}(\Theta-\theta_-)dxd\tau\nonumber\\
&&=a\int_0^t(\ln\Theta)_x(-st,t)sgn_{\eta}(\Theta-\theta_-)(-st,t)d\tau.
\end{eqnarray}

Again, both side of (\ref{1.8})$_1$ multiply by
$sgn_{\eta}(\Theta-\theta_+)$ and integrate in
$(-s\tau,+\infty)\times (0,t)$ we can get
\begin{eqnarray*}
\int_0^t\left(\int_{-s\tau}^{+\infty}I_{\eta}(\Theta-\theta_+)dx\right)_{\tau}
d\tau&=&-a\int_0^t(\ln\Theta)_x(-st,t)sgn_{\eta}(\Theta-\theta_+)(-st,t)d\tau\\
&&\quad-a\int_0^t\int_{-s\tau}^{+\infty}(\ln\Theta)_x^2sgn'_{\eta}(\Theta-\theta_+)dxd\tau.
\end{eqnarray*}
When $\eta\to 0$ and use (\ref{2.11}) , (\ref{2.15}) and Lemma
\ref{lem2.1} we can get
\begin{eqnarray}\label{2.47}
&&\|\Theta-\theta_+\|_{L^1(\Omega_2)}
+a\int_0^t\int_{-s\tau}^{+\infty}(\ln\Theta)_x^2sgn'_{\eta}(\Theta-\theta_+)dxd\tau\nonumber\\
&&=-a\int_0^t(\ln\Theta)_x(-st,t)sgn_{\eta}(\Theta-\theta_-)(-st,t)d\tau+\|\Theta_0-\theta_+\|_{L^1(\Omega_2)}.
\end{eqnarray}
Similar as (\ref{2.34}), when we integrate (\ref{2.30}) in
$\R_+\times (0,t)$ and combine with (\ref{2.11})
 we can get there exist constant $C>0$ independent of
$\alpha$ such that
 \begin{equation}\label{2.48}
\|(\ln\Theta)_{xx}\|^2+a\int_0^t\|\partial_x^3(\ln\Theta)\|^2d\tau+\int_0^t(\ln\Theta)_{xx}^2(0,\tau)d\tau\leq
Ca^{-1}\alpha^2+C\alpha^3+Ca\alpha^3+Ca^{-1}\alpha.
 \end{equation}
Since $a=\kappa p_+(\gamma-1)/(\gamma R^2)$ and $\kappa\to 0$,we can
choose $\kappa=\alpha^{-1/2}<1$, use (\ref{2.11}) and (\ref{2.48})
such that (\ref{2.46}) and (\ref{2.47}) are meant
$$\|\Theta-\theta_-\|_{L^1(\Omega_1)}+\|\Theta-\theta_+\|_{L^1(\Omega_2)}\leq Ct(a\alpha+(a\alpha)^{3/4}+a^{5/4}\alpha)+C\alpha^{-1}\leq C\kappa^{3/8}(t+1),$$
so we get $\|(V-\overline{V},\Theta-\overline{\Theta})\|_{L^p}\to 0$
as $\kappa\to 0$ with any $t\in [0,T]$.

It remains to estimate $\|U-\overline{U}\|_{L^p}$. To do so, both
side of (\ref{2.29}) multiply by $sgn_{\eta}((\ln\Theta)_x)$ then
integrate in $\R_+\times(0,t)$ we can get
\begin{eqnarray*}&&\int_0^t\left(\int_{\R_+}I_{\eta}((\ln\Theta)_x)dx\right)_{\tau}d\tau+s\int_0^t
I_{\eta}((\ln\Theta)_x)(0,\tau)d\tau+a\int_0^t\int_{\R_+}\frac{(\ln\Theta)^2_{xx}}{\Theta}sgn'_{\eta}((\ln\Theta)_x)dxd\tau\\
&&=-a\int_0^t\theta_-^{-1}(\ln\Theta)_{xx}(0,\tau)sgn_{\eta}((\ln\Theta)_x)(0,\tau)d\tau.\end{eqnarray*}
Again let $\eta\to 0$, $\kappa=\alpha^{-1/2}<1$ we can get from
(\ref{2.11}), (\ref{2.48}) and Lemma \ref{lem2.1} that there exist
constant $C>0$ independent of $\alpha$ such that
\begin{eqnarray*}&&\int_{\R_+}|(\ln\Theta)_x|dx+a\int_0^t\int_{\R_+}\frac{(\ln\Theta)^2_{xx}}{\Theta}sgn'_{\eta}((\ln\Theta)_x)dxd\tau\\
&&\leq
Ct^{1/2}\left(\int_0^t(\ln\Theta)^2_{xx}(0,\tau)d\tau\right)^{1/2}+Ct^{1/2}\left(\int_0^t(\ln\Theta)^2_{x}(0,\tau)d\tau\right)^{1/2}
+\|(\ln\Theta_0)_x\|_{L^1}\\
&&\leq
Ct^{1/2}(a^{-1/2}\alpha+\alpha^{3/2}+a^{1/2}\alpha^{3/2}+\alpha^{1/2})+C\leq
C(1+t)^{1/2}\alpha^2.\end{eqnarray*}

Use the definition of $U$ in (\ref{1.7}) and combine with
(\ref{2.11}),(\ref{2.48}) and $\kappa=\alpha^{-1/2}<1$ we know that
\begin{eqnarray*}
&&\|U-\overline{U}\|^p_{L^p}\leq
C\kappa^p\|(\ln\Theta)_x\|_{L^1}\|(\ln\Theta)_x\|^{(p-1)/2}\|(\ln\Theta)_{xx}\|^{(p-1)/2}\\
&&\leq
C(\alpha^{5(p-1)/4}+\alpha^{p-1}+\alpha^{(3-p)/2})\alpha^{-2p}\alpha^2(1+t)^{1/2}\\
&&\leq
C\left(\alpha^{-3(p-1)/4}+\alpha^{1-p}+\alpha^{-5(p-1)/2}\right)(1+t)^{1/2}.
\end{eqnarray*}
Remind that $\alpha=\kappa^{-1/2}$, so we can get
$$\lim_{\kappa\to 0}\|U-\overline{U}\|_{L^p}=0.$$
The proof of Lemma \ref{lem2.4} is therefore complete, which also
means $(V,U,\Theta)$ is viscous contact discontinuity.$\Box$

\section{Proof of Theorem \ref{thm1.1}}

   Under the preparations in last section, the main task here is to finish (\ref{2.1}).  This part we also do some preparations. we must use the results
\begin{eqnarray}\label{3.1}
&&|V_x|\leq C|\Theta_x|,\nonumber\\
&&|\Theta_x|^2\leq
C\|(\ln\Theta)_x\|\|(\ln\Theta)_{xx}\|,\nonumber\\
&&|U_x|\leq C |(\ln\Theta)_{xx}|,\nonumber\\
 &&|U_x|^2\leq
C\|(\ln\Theta)_{xx}\|\|(\ln\Theta)_{xxx}\|,
\end{eqnarray}
 which come from
(\ref{1.7})--(\ref{1.9}) . Also we set $C(\delta_0)$ stands for
small constants about  $\delta_0$, $\|\varphi_0,\psi_0,\zeta_0\|$ is
asked suitably small, $C_v=\frac{R}{\gamma-1}$ and
$$\epsilon_1\ll \epsilon_3\ll \epsilon_2.$$

Now , let's finish (\ref{2.1}) which is very important for our
 proof of Theorem \ref{thm1.1}.

\begin{lem}\label{lem3.1}
If  $\epsilon_1>0$  and  $C(\delta_0)>0$ are small constant about
$\delta_0$,  we can get
\begin{eqnarray*}
&&\int_{\R_+}\left(R\theta\Phi\left(\frac{v}{V}\right)+\frac{1}{2}\psi^2+C_v\theta\Phi\left(\frac{\theta}{\Theta}\right)\right)dx+\int_0^t\
\left\|\left(\frac{\psi_x}{\sqrt{v\theta}},\frac{\zeta_x}{\theta\sqrt{v}}\right)\right\|^2\
d\tau\\
&&\leq
C\epsilon_1^{-1}\int_0^t\int_0^{+\infty}\Theta_x^2(\varphi^2+\zeta^2)\
dxd\tau+C(\delta_0)+C\left\{\epsilon_1\int_0^t\
\left(\|\varphi_x\|^2+\psi_x^2(0,\tau)\right)\
d\tau+\|(\varphi_0,\psi_0,\zeta_0)\|^2\right\}.
\end{eqnarray*}
\end{lem}
\pf Set
$$\Phi(z)=z-\ln z-1,$$
$$\Psi(z)=z^{-1}+\ln z-1,$$ where $\Phi'(1)=\Phi(1)=0$ is a strictly
convex function around $z=1$. Similar to the proof in \cite{HMS},
 we deduce from (\ref{1.12}) that
\begin{eqnarray}\label{3.2}
&&\left(\frac{\psi^2}{2}+R\Theta\Phi\left(\frac{v}{V}\right)+C_v\Theta\Phi\left(\frac{\theta}{\Theta}\right)\right)_t-
s\big(\frac{\psi^2}{2}+R\Theta\Phi\left(\frac{v}{V}\right)+C_v\Theta\Phi\left(\frac{\theta}{\Theta}\right)\big)_x\nonumber\\
&&\quad+\mu\frac{\Theta\psi_x^2}{v\theta}+\kappa\frac{\Theta\zeta_x^2}{v\theta^2}+H_x+Q
=\mu\left(\frac{\psi\psi_x}{v}\right)_x-F\psi-\frac{\zeta
G}{\theta},
\end{eqnarray}
where
$$H=R\frac{\zeta\psi}{v}-R\frac{\Theta\varphi\psi}{vV}+\mu\frac{U_x\varphi\psi}{vV}-\kappa\frac{\zeta\zeta_x}{v\theta}+\kappa\frac{\Theta_x\varphi\zeta}{v\theta V},$$
and
\begin{eqnarray*}
Q&=&p_+\Phi\left(\frac{V}{v}\right)U_x+\frac{p_+}{\gamma-1}\Phi\left(\frac{\Theta}{\theta}\right)U_x-\frac{\zeta}{\theta}(p_+-p)U_x-\mu\frac{U_x\varphi\psi_x}{vV}\\
&&\quad-\kappa
\frac{\Theta_x}{v\theta^2}\zeta\zeta_x-\kappa\frac{\Theta\Theta_x}{v\theta^2V}\varphi\zeta_x-2\mu\frac{U_x}{v\theta}\zeta\psi_x
+\kappa\frac{\Theta_x^2}{v\theta^2V}\varphi\zeta+\mu\frac{U_x^2}{v\theta
V}\varphi\zeta\\
&=:&\sum_{i=1}^9Q_i.
\end{eqnarray*}
Note that $p=R\theta/v$, $p_+=R\Theta/V$ and (\ref{1.7}), use
integrate by part and Cauchy-Schwarz inequality can get

\begin{eqnarray}\label{3.3}
Q_1+Q_2&=&Ra\left(\Phi\left(\frac{V}{v}\right)(\ln\Theta)_x\right)_x+\frac{Ra}{\gamma-1}\left(\Phi\left(\frac{\Theta}{\theta}\right)(\ln\Theta)_x\right)_x\nonumber\\
&&\quad-aR(\ln\Theta)_x\left(\frac{V\varphi_x\varphi-V_x\varphi^2}{Vv^2}\right)\nonumber\\
&&\quad-a\frac{p_+}{\gamma-1}(\ln\Theta)_x\left(\frac{\Theta\zeta_x\zeta-\Theta_x\zeta^2}{\Theta\theta^2}\right)\nonumber\\
&&\geq\left(p_+\Phi\left(\frac{V}{v}\right)U+\frac{p_+}{\gamma-1}\Phi\left(\frac{\Theta}{\theta}\right)U\right)_x\nonumber\\
&&\quad-\epsilon(\zeta_x^2+\varphi_x^2)-C\epsilon^{-1}\Theta_x^2(\zeta^2+\varphi^2).
\end{eqnarray}
Similarly, using $p-p_+=\frac{R\zeta-p_+\varphi}{v}$, we can get
\begin{equation}\label{3.4}
Q_3\geq\frac{R\zeta-p_+\varphi}{v}(\frac{\zeta}{\theta}U_x)\geq\left(\frac{R\zeta^2U}{v\theta}-\frac{p_+\zeta\varphi
U}{\theta
v}\right)_x-\epsilon(\zeta_x^2+\varphi_x^2)-C\epsilon^{-1}\Theta_x^2(\zeta^2+\varphi^2).
\end{equation}
And
\begin{eqnarray}\label{3.5}
(Q_4+Q_7)+(Q_5+Q_6+Q_8)+Q_9&\geq&
-C\epsilon^{-1}(\ln\Theta)_{xx}^2-\epsilon\psi_x^2\nonumber\\
&&-\epsilon\zeta_x^2-C\epsilon^{-1}\Theta_x^2(\zeta^2+\varphi^2)\nonumber\\
&&-C\epsilon^{-1}|(\ln\Theta)_{xx}|^2(\zeta^2+\varphi^2).
\end{eqnarray}
At the end we use the definition of $F$ and $G$ in (\ref{1.10}) then
combine with the general inequality skills as above to get
\begin{eqnarray}\label{3.6}
-F\psi-G\frac{\zeta}{\theta}&=&-\frac{\kappa a(\gamma-1)-\mu
p_+\gamma}{R\gamma}\left(\frac{(\ln\Theta)_{xx}}{\Theta}\right)_x\psi\nonumber\\
&&\quad+\frac{\mu p_+}{R\Theta}\left(\frac{\kappa
(\gamma-1)}{R\gamma}(\ln\Theta)_{xx}\right)^2\frac{\zeta}{\theta}\nonumber\\
&\leq&-\frac{\kappa a(\gamma-1)-\mu
p_+\gamma}{R\gamma}\left(\frac{(\ln\Theta)_{xx}}{\Theta}\psi\right)_x+\frac{\kappa
a(\gamma-1)-\mu
p_+\gamma}{R\gamma}\frac{(\ln\Theta)_{xx}}{\Theta}\psi_x\nonumber\\
&&\quad+\frac{\mu p_+}{R\Theta}\left(\frac{\kappa
(\gamma-1)}{R\gamma}(\ln\Theta)_{xx}\right)^2\frac{\zeta}{\theta}\nonumber\\
&\leq&-\frac{\kappa a(\gamma-1)-\mu
p_+\gamma}{R\gamma}\left(\frac{(\ln\Theta)_{xx}}{\Theta}\psi\right)_x+\epsilon\psi_x^2+C\epsilon^{-1}(\ln\Theta)_{xx}^2.
\end{eqnarray}
Integrating (\ref{3.3})$-$(\ref{3.6}) in $\R\times(0,t)$ ,  using
(\ref{2.11}), (\ref{2.15}) and the boundary condition about
$(\varphi,\psi,\zeta)$ of (\ref{1.12}) to estimate the terms
$\mu\left(\frac{\psi\psi_x}{v}\right)_x$,
$\left(\frac{(\ln\Theta)_{xx}\psi}{\Theta}\right)_x$ and $H_x$, in
the end combine with Cauchy-Schwarz inequality  we now   that for a
small $\epsilon>0$ which is about $\delta_0$,
$C_v=\frac{R}{\gamma-1}$, we have
\begin{eqnarray}\label{3.7}
&&\int_{\R_+}\left(R\theta\Phi\left(\frac{v}{V}\right)+\frac{1}{2}\psi^2+C_v\theta\Phi\left(\frac{\theta}{\Theta}\right)\right)dx+\int_0^t\
\left\|\left(\psi_x/(\sqrt{v\theta}),\zeta_x/(\theta\sqrt{v})\right)\right\|^2\
d\tau\nonumber\\
&&\leq
C\epsilon^{-1}\left\{\int_0^t\int_0^{+\infty}\Theta_x^2(\varphi^2+\zeta^2)\
dxd\tau+\|\Theta_{0x}\|^2\right\}+C\left\{\epsilon\int_0^t\ \|\varphi_x\|^2\ d\tau+\|(\varphi_0,\psi_0,\zeta_0)\|^2\right\}\nonumber\\
&&\qquad+C\int_0^t\psi^2(0,\tau)d\tau+\epsilon\int_0^t\psi_x^2(0,\tau)
d\tau+C\int_0^t(\ln\Theta)_{xx}^2(0,\tau)d\tau+C(\delta_0).
\end{eqnarray}

Using the definition about $\psi(0,t)$ in (\ref{1.12}), then combine
with (\ref{1.9})$_5$ and (\ref{2.11}), (\ref{2.15}) we can get
\begin{eqnarray}
\int_0^t\psi^2(0,\tau)d\tau+\int_0^t(\ln\Theta)_{xx}^2(0,\tau)d\tau\leq
C(\delta_0),\label{3.8}
\end{eqnarray}
Insert (\ref{3.8}) into (\ref{3.7})  we finish this lemma.$\Box$

\begin{lem}\label{lem3.2} If $C(\delta_0)>0$ is a small constant about
$\delta_0$
\begin{eqnarray*}&&\int_0^t\int_{\R_+}\Theta_x^2(\varphi^2+\zeta^2)dxd\tau\leq
C(\delta_0)\int_0^t\|(\varphi_x,\zeta_x)\|^2d\tau.
\end{eqnarray*}
\end{lem}
\pf Because if $x>0$
\begin{eqnarray*}&&\frac{\varphi^2}{x+1}=\int_0^x\big(\frac{2\varphi\varphi_x}{x+1}-\frac{\varphi^2}{(x+1)^2}\big)dx\\
&&=\int_0^x\left(\varphi_x^2-(\varphi_x-\frac{\varphi}{x+1})^2\right)dx\leq\int_0^x\varphi_x^2dx\leq
\|\varphi_x\|^2.\end{eqnarray*} similar as above we can
get\begin{eqnarray*}&&\frac{\zeta^2}{x+1}\leq\int_0^x\zeta_x^2dx\leq
\|\zeta_x\|^2.\end{eqnarray*}

As to
\begin{eqnarray*}\int_0^t\int_0^{+\infty}\Theta_x^2(\varphi^2+\zeta^2)dxd\tau&\leq&
\int_0^t\int_0^{+\infty}\Theta_x^2(x+1)\frac{(\varphi^2+\zeta^2)}{1+x}dxd\tau\\
&&\leq
\int_0^t\big(\int_0^{+\infty}\Theta_x^2(1+x)dx\big)\|(\varphi_x,\zeta_x)\|^2d\tau,
\end{eqnarray*}
 use (\ref{2.11})and (\ref{2.18}) we can get \begin{eqnarray*}&&\int_0^t\int_{\R_+}\Theta_x^2(\varphi^2+\zeta^2)dxd\tau\leq C(\delta_0)\int_0^t\|(\varphi_x,\zeta_x)\|^2d\tau,
\end{eqnarray*} and we finish this lemma.$\Box$

\begin{lem}\label{lem3.3} If a constant $\epsilon_2>0$ , $C(\delta_0)>0$ is a small constant about
$\delta_0$, we can get
\begin{eqnarray*}
&&\|(\varphi,\psi,\zeta)\|^2+\|(\psi_x,\zeta_x)\|^2+\int_0^t(\psi^2_x(0,\tau)+\zeta^2_x(0,\tau))d\tau+\int_0^t\|(\psi_{xx},\zeta_{xx})\|^2d\tau\\
&&\leq
C\left(\|(\psi_{0x},\zeta_{0x})\|^2+\epsilon_2^{-1}\|(\varphi_0,\psi_0,\zeta_0)\|^2\right)+C\epsilon_2^{-1}\int_0^t\|\varphi_x\|^2d\tau+C(\delta_0).
\end{eqnarray*}.
\end{lem}

\pf First to get the estimate of $\|\psi_x(t)\|$ ,multiply both side
of  (\ref{1.12})$_2$ to $\psi_{xx}$ to get
\begin{eqnarray*}
&&\left(\frac{\psi_x^2}{2}\right)_t+s\left(\frac{\psi_x^2}{2}\right)_x+\mu\frac{\psi_{xx}^2}{v}
=\mu\frac{\psi_x v_x}{v^2}\psi_{xx}+\mu\left(\frac{U_x\varphi}{v
V}\right)_x\psi_{xx}\\
&&\quad-R\left(\frac{\Theta\varphi}{v
V}\right)_x\psi_{xx}+R\left(\frac{\zeta}{v}\right)_x\psi_{xx}+F\psi_{xx}+(\psi_t\psi_x)_x:=\sum_{i=1}^6I_i.
\end{eqnarray*}
use last inequality integrate in $\R_+\times (0,t)$ ($s=-u_b/v_-<0$)
to get
\begin{eqnarray}
&&\|\psi_x(t)\|^2+\int_0^t\psi_x^2(0,\tau)d\tau+\int_0^t\|\psi_{xx}(\tau)\|^2d\tau\nonumber\\
&&\quad\leq
C\|\psi_{0x}\|^2+C\sum_{i=1}^6\left|\int_0^t\int_0^\infty
I_idxd\tau\right|.\label{3.9}
\end{eqnarray}

Now deal with   $\iint|I_i|dxd\tau$ in the right side of
(\ref{3.9}). Using $\epsilon$ small and $v=\varphi+V$ ,
$R\Theta/V=p_+$ and (\ref{2.13}), (\ref{2.14}), (\ref{3.1})to get
\begin{eqnarray}
&&\int_0^t \int_0^{+\infty}|I_1| dx d\tau\leq C\int_0^t
\int_0^{+\infty}|V_x|
|\psi_{x}||\psi_{xx}|dxd\tau+C\int_0^t\int_0^{+\infty}
|\varphi_x||\psi_{x}||\psi_{xx}|dx d\tau\nonumber\\
&&\quad\leq C\int_0^t\|V_x\|\|\psi_x\|_{L^\infty}\|\psi_{xx}\|
d\tau+C\int_0^t \|\psi_x\|_{L^\infty}\|\varphi_x\|\|\psi_{xx}\|
d\tau\nonumber\\
&&\quad\leq \epsilon\int_0^t\|\psi_{xx}\|^2d\tau+
C\epsilon^{-1}\int_0^t\|\psi_x\|^2 \|V_x\|^4d\tau+C\int_0^t
\|\psi_x\|^{1/2}\|\varphi_x\|\|\psi_{xx}\|^{3/2}
d\tau\nonumber\\
&&\quad\leq C\epsilon\int_0^t\|\psi_{xx}\|^2d\tau+
C(\delta_0)\int_0^t\|\psi_x\|^2
d\tau\nonumber\\
&&\quad\quad+C\epsilon^{-1}\sup_{t}\|\varphi_x\|^4\int_0^t
\|\psi_x\|^2 d\tau .\label{3.10}
\end{eqnarray}
Next we use the definition of $(V,U,\Theta)$ in
(\ref{1.7}),(\ref{1.9}), Cauchy-Schwarz inequality and (\ref{2.15}),
(\ref{2.16}), (\ref{3.1}) to get
\begin{eqnarray}
&&\int_0^t\int_0^\infty|I_2|dxd\tau\nonumber\\
&&\quad\leq
C\int_0^t\int_0^\infty\left(|U_{xx}||\varphi|+|U_x||\varphi_x|+|U_x||V_x||\varphi|+|U_x||\varphi||\varphi_x|\right)|\psi_{xx}|dxd\tau\nonumber\\
&&\quad\leq\epsilon\int_0^t\|\psi_{xx}\|^2d\tau+\frac{C}{\epsilon}\int_0^t\|\varphi\|_{L^\infty}^2\|U_{xx}\|^2d\tau
+\frac{C}{\epsilon}\int_0^t\|U_x\|_{L^\infty}^2\|\varphi_x\|^2d\tau\nonumber\\
&&\qquad+\frac{C}{\epsilon}\int_0^t\|\varphi\|_{L^\infty}^2\|V_x\|^2\|U_x\|_{L^\infty}^2d\tau
+\frac{C}{\epsilon}\int_0^t\|\varphi\|_{L^\infty}^2\|U_x\|_{L^\infty}^2\|\varphi_x\|^2d\tau\nonumber\\
&&\quad\leq\epsilon\int_0^t\|\psi_{xx}\|^2d\tau+C(\delta_0)+C(\delta_0)\int_0^t\|\varphi_x\|^2d\tau.\label{3.11}
\end{eqnarray}
The same as (\ref{3.10}) and (\ref{3.11}), we use Lemma
\ref{lem3.2}, the definition of $F$ in(\ref{2.2}) and
(\ref{2.11}),(\ref{2.13})--(\ref{2.15}),(\ref{3.1})we can get the
estimates about $I_3$ to $I_5$ as following.
\begin{eqnarray}
&&\int_0^t\int_0^\infty(|I_3|+|I_4|+|I_5|)dxd\tau\nonumber\\
&&\leq
C\int_0^t\int_0^\infty\left(|\Theta_{x}||\varphi|+|\Theta||\varphi_x|+|\Theta||V_x||\varphi|
+|\Theta||\varphi||\varphi_x|\right)|\psi_{xx}|dxd\tau\nonumber\\
&&\quad+C\int_0^t\int_0^\infty\left(|\zeta_x|+|\zeta||V_x|+|\zeta||\varphi_x|\right)|\psi_{xx}|dxd\tau\nonumber\\
&&\quad+\epsilon\int_0^t\|\psi_{xx}\|^2d\tau+\frac{C}{\epsilon}\int_0^t\|F\|^2d\tau\nonumber\\
&&\leq\epsilon\int_0^t\|\psi_{xx}\|^2d\tau+C\epsilon^{-1}\int_0^t\|\varphi_x\|^2d\tau+\frac{C}{\epsilon}\int_0^t\|\varphi_x\|^2d\tau+\frac{C}{\epsilon}\int_0^t\int_0^\infty
V_x^2\varphi^2dx
d\tau\nonumber\\
&&\quad+\epsilon\int_0^t\|\psi_{xx}\|^2d\tau+\frac{C}{\epsilon}\int_0^t\int_0^\infty\left(\zeta_x^2+V_x^2\zeta^2\right)dxd\tau+\frac{C}{\epsilon}\sup_t\|(\varphi,\zeta)\|\|(\varphi_x,\zeta_x)\|\int_0^t\|\varphi_x\|^2d\tau\nonumber\\
&&\quad+\epsilon\int_0^t\|\psi_{xx}\|^2d\tau+C(\delta_0).\label{3.12}
\end{eqnarray}

Because Lemma \ref{lem3.1} and Lemma \ref{lem3.2}, we know
$\|(\varphi,\zeta)\|$ is suitably small when $C(\delta_0)$ and
$\|(\varphi_0,\psi_0,\zeta_0)\|$ small. So there exist a small
constant $\delta$ about $\|(\varphi_0,\psi_0,\zeta_0)\|$ and
$\delta_0$  such that
$$\frac{C}{\epsilon}\sup_t\|(\varphi,\zeta)\|\|(\varphi_x,\zeta_x)\|\int_0^t\|\varphi_x\|^2d\tau\leq C\delta\int_0^t\|\varphi_x\|^2d\tau+C\delta\int_0^t\psi_x^2(0,\tau)d\tau,$$
here $\|\varphi_x\|^2+\int_0^t\|\varphi_x\|^2d\tau\leq C$ can be
established in Lemma \ref{lem3.4}.

 Therefore
\begin{eqnarray*}
&&\int_0^t\int_0^\infty(|I_3|+|I_4|+|I_5|)dxd\tau\\
&&\leq
\epsilon\int_0^t\|\psi_{xx}\|^2d\tau+\frac{C}{\epsilon}\int_0^t\int_0^\infty\
(\zeta_x^2+\varphi_x^2)dxd\tau
+C\delta\int_0^t\psi_x^2(0,\tau)d\tau+C(\delta_0).\end{eqnarray*}

At last we integrate by part to the term about $I_6$ to get
\begin{eqnarray}
&&\left|\int_0^t\int_0^\infty
I_6dxd\tau\right|=\left|\int_0^t(\psi_t\psi_x)(0,\tau)d\tau\right|\leq
\frac{C}{C^{1/2}(\delta_0)}\int_0^t\psi_x^2(0,\tau)d\tau+C^{1/2}(\delta_0)\int_0^t\psi_\tau^2(0,\tau)d\tau\nonumber\\
&&\leq \frac{C}{\epsilon}C^{-1/2}(\delta_0)\int_0^t \|\psi_x\|^2\
d\tau+\epsilon\int_0^t \|\psi_{xx}\|^2\
d\tau+C^{1/2}(\delta_0)\int_0^t\psi_\tau^2(0,\tau)d\tau.\label{3.13}
\end{eqnarray}
Using the definition of $U$ in (\ref{1.7}), $\psi=u-U$ and
(\ref{2.29}) to get
\begin{eqnarray}
\psi_t(0,t)&=&-\frac{k(\gamma-1)}{\gamma
R}(\ln\Theta)_{xt}(0,t)\nonumber\\
&&=-s\frac{k(\gamma-1)}{\gamma
R}(\ln\Theta)_{xx}(0,t)-a\frac{k(\gamma-1)}{\gamma
R}\partial_x\left(\frac{(\ln\Theta)_{xx}}{\Theta}\right)(0,t).
\label{3.14}
\end{eqnarray}

Combine with (\ref{2.15}) (\ref{2.17}) and $|\Theta_x(0,t)|\leq C$
we get
\begin{equation}
\int_0^t\|(\ln\Theta)_{x\tau}\|_{L^\infty}^2(0,\tau)d\tau\leq
C.\label{3.15}
\end{equation}
So combine with (\ref{3.13}) and (\ref{3.15}) we get
\begin{eqnarray}
&&\left|\int_0^t\int_0^\infty I_6dxd\tau\right|\nonumber\\
&&\leq
 \frac{C}{\epsilon}\int_0^t \|\psi_x\|^2\
d\tau+\epsilon\int_0^t \|\psi_{xx}\|^2\
d\tau+C(\delta_0).\label{3.16}
\end{eqnarray}

In all there exist a small $\delta>0$
\begin{eqnarray}\label{3.17}
&&\int_0^t\ \int_0^{+\infty}\ \sum_{i=1}^6|I_i|\ dx\ d\tau\nonumber\\
&&\leq C\int_0^t\
\left(\epsilon\|\psi_{xx}\|^2+\delta\psi_x^2(0,\tau)\right)\ d\tau+CN^4(t)\epsilon^{-1}\int_0^t\|\psi_x\|^2d\tau\nonumber\\
&&\quad+C\epsilon^{-1}\int_0^t\ \|(\varphi_x,\psi_x,\zeta_x)\|^2\
d\tau+C(\delta_0).
\end{eqnarray}

So (\ref{3.9}) can be change to
\begin{eqnarray}
&&\|\psi_x(t)\|^2+\int_0^t\psi_x^2(0,\tau)d\tau+\int_0^t\|\psi_{xx}(\tau)\|^2d\tau\nonumber\\
&&\leq
C(\delta_0)+C\epsilon^{-1}\int_0^t\|(\varphi_x,\psi_x,\zeta_x)\|^2d\tau+CN^4(t)\epsilon^{-1}\int_0^t\|\psi_x\|^2d\tau+C\|\psi_{0x}\|^2.
\label{3.18}
\end{eqnarray}

The estimate about  $\|\zeta_x\|$ is similar to $\|\psi_x\|$, use
(\ref{1.12})$_3$ multiply $\zeta_{xx}$ then integrate in
$Q_t=\R_+\times(0,t)$ to get
\begin{eqnarray}
&&\|\zeta_x\|^2+\int_0^t\zeta_x^2(0,\tau)d\tau+\int_0^t
\|\zeta_{xx}\|^2\ d\tau\nonumber\\
&&\leq C \|\zeta_{0x}\|^2+
C\epsilon^{-1}\int_0^t\int_0^\infty\left(\psi_x^2+\zeta\psi_x^2+\zeta^2U_x^2+U_x^2\varphi^2\right)dxd\tau\nonumber\\
&&
\quad+C\int_0^t\int_0^{+\infty}|\zeta_x|(|\varphi_x|+|V_x|)|\zeta_{xx}|dxd\tau+C\epsilon^{-1}\int_0^t\int_0^\infty\left|\left(\frac{\Theta_x\varphi}{vV}\right)_x\right|^2dxd\tau
\nonumber\\
&&\quad+C\epsilon^{-1}\int_0^t\int_0^{+\infty}(U_x^4+\psi_x^4)dxd\tau+C\epsilon^{-1}\int_0^t\|G\|^2d\tau\nonumber\\
&&=:C\|\zeta_{0x}\|^2+\sum_{i=1}^5J_i.\label{3.19}
\end{eqnarray}
Use the same method as (\ref{3.10})--(\ref{3.13})
$$
J_1\leq
C\epsilon^{-1}\int_0^t\|\psi_x\|^2d\tau+C\epsilon^{-1}N^2(t)\int_0^t\|U_x\|^2d\tau
\leq C\epsilon^{-1}\int_0^t\|\psi_x\|^2d\tau+C(\delta_0).
$$
Again use the same method as(\ref{3.10})--(\ref{3.13})
\begin{eqnarray*}
J_2&\leq& C\int_0^t\|\zeta_x\|_{L^\infty}\|\varphi_x\|\|\zeta_{xx}\|d\tau+C\int_0^t\|V_x\|\|\zeta_x\|_{L^\infty}\|\zeta_{xx}\|d\tau\\
&\leq&C\int_0^t\|\zeta_x\|^{1/2}\|\zeta_{xx}\|^{3/2}\|\varphi_x\|d\tau
+\epsilon\int_0^t\|\zeta_{xx}\|^2d\tau+C(\delta_0)\int_0^t\|\zeta_x\|^2d\tau\\
&\leq&2\epsilon\int_0^t\|\zeta_{xx}\|^2d\tau+C(\delta_0)\int_0^t\|\zeta_x\|^2d\tau
+C\epsilon^{-1}\sup_t\|\varphi_x\|^4\int_0^t\|\zeta_x\|^2d\tau.
\end{eqnarray*}

Because
\begin{eqnarray*}&&\left|\left(\frac{\Theta_x\varphi}{vV}\right)_x\right|^2\\
&&=|\frac{\Theta_{xx}\varphi}{vV}+\frac{\Theta_x\varphi_x}{vV}+\frac{\Theta_x\varphi}{vV}(-\frac{V_x+\varphi_x}{v^2}-\frac{V_x}{V^2})|^2\\
&&\leq
C\Theta_{xx}^2\varphi^2+C\Theta_x^2\varphi_x^2+C\Theta_x^2V_x^2\varphi^2+C\Theta_x^2\varphi^2\varphi_x^2,
\end{eqnarray*}
combine with $R\Theta/V=p_+$, use the same method as
(\ref{3.10})--(\ref{3.13}) to get
\begin{eqnarray*}
J_3&\leq&
C\epsilon^{-1}\int_0^t\|\varphi\|_{L^\infty}^2\|\Theta_{xx}\|^2d\tau
+C\epsilon^{-1}\int_0^t\|\Theta_x\|_{L^\infty}^2\|\varphi_x\|^2d\tau\\
&&\quad+C\epsilon^{-1}\int_0^t\ \int_0^{+\infty}\ \Theta_x^2V_x^2\varphi^2\ dx\ d\tau\nonumber\\
&\leq&C(\delta_0)\int_0^t\|\varphi_{x}\|^2d\tau+C(\delta_0) .
\end{eqnarray*}
Use the definition $U$ and similar as (\ref{3.10}) (\ref{3.11})
that we combine with Lemma \ref{lem2.3} to get
\begin{eqnarray*} J_4&\leq&
C(\delta_0)+C\epsilon^{-1}\int_0^t\|\psi_{x}\|_{L^\infty}^2\|\psi_x\|^2d\tau\\
&\leq& C(\delta_0)+C\epsilon^{-1}\int_0^t\|\psi_x\|^3\|\psi_{xx}\|d\tau\\
&\leq&
C(\delta_0)+C\int_0^t\left(\epsilon^{-2}\|\psi_x\|^2\|\psi_x\|^4+\epsilon^{2}\|\psi_{xx}\|^2\right)d\tau.
\end{eqnarray*}
Use the definition $G$ in (\ref{2.2}) combine with Lemma
\ref{lem2.3}
\begin{eqnarray*}
J_5=C\epsilon^{-1}\int_0^t\|G\|^2d\tau\leq C(\delta_0).
\end{eqnarray*}
Use the results from $J_1$ to $J_5$, the inequality (\ref{3.19}) can
be change to
\begin{eqnarray}
&&\|\zeta_x\|^2+\int_0^t\zeta_x^2(0,\tau)d\tau+\int_0^t\ \|\zeta_{xx}\|^2\ d\tau\nonumber\\
&&\leq
C\|\zeta_{0x}\|^2+C(\epsilon^{-3}+N^4(t))\int_0^t\|(\psi_x,\zeta_x)\|^2d\tau+C(\delta_0)\int_0^t\|\varphi_x\|^2d\tau+C(\delta_0)\nonumber\\
&&\quad+C\epsilon\int_0^t\|\psi_{xx}\|^2d\tau.\label{3.20}
\end{eqnarray}

 In fact when combine  with Lemma \ref{lem3.1}--\ref{lem3.2}, (\ref{3.18}) and (\ref{3.20}), it is easy to get
\begin{eqnarray*}
&&\|(\varphi,\psi,\zeta)\|^2+\|(\psi_x,\zeta_x)\|^2+\int_0^t(\psi^2_x(0,\tau)+\zeta^2_x(0,\tau))d\tau+\int_0^t\|(\psi_{xx},\zeta_{xx})\|^2d\tau\\
&&\leq
C\left(\|(\psi_{0x},\zeta_{0x})\|^2+\epsilon^{-3}\|(\varphi_0,\psi_0,\zeta_0)\|^2\right)+C\epsilon^{-1}\int_0^t\|\varphi_x\|^2d\tau+C(\delta_0).
\end{eqnarray*}
$\Box$

\begin{lem}\label{lem3.4} For a small $\epsilon_3>0$  and $C(\delta_0)>0$ is a small constant about
$\delta_0$  , we can get
\begin{eqnarray}
\|\varphi_x\|^2+\int_0^t\|\varphi_x\|^2 d\tau&\leq&
C\|\varphi_{0x}\|^2+C\epsilon_3^{-1}\|(\varphi_0,\psi_0,\zeta_0)\|^2
+\int_0^t\epsilon_3\|\psi_{xx}\|^2d\tau\nonumber\\
&&\quad+\int_0^tC\epsilon_3^{-1}\|(\psi_x,\zeta_x)\|^2d\tau+C(\delta_0)
.\label{3.21}
\end{eqnarray}
\end{lem}
\pf Set $\bar{v}=\frac{v}{V}$ take it into (\ref{1.12})$_1$,
(\ref{1.12})$_2$  ($p=R \theta/v$) to get
$$
\psi_t-s\psi_x+p_x=\mu\left
(\frac{\bar{v}_x}{\bar{v}}\right)_t-s\mu\left
(\frac{\bar{v}_x}{\bar{v}}\right)_x-F,
$$
Both sides of last equation multiply $\bar v_x/\bar v$ to get
\begin{eqnarray}
&&\left(\frac{\mu}{2}\left (\frac{\bar v_x}{\bar
v}\right)^2-\psi\frac{\bar v_x}{\bar v}\right)_t
+\frac{R\theta}{v}\left (\frac{\bar v_x}{\bar
v}\right)^2+\left(\psi\frac{\bar v_t}{\bar v}\right)_x-s\mu\left
(\frac{\bar v_x}{\bar v}\right)_x\frac{\bar v_x}{\bar v}\nonumber\\
&&\quad=\frac{\psi_x^2}{v}+U_x\left
(\frac{1}{v}-\frac{1}{V}\right)\psi_x+\frac{R\zeta_x}{v}\frac{\bar
v_x}{\bar v}-\frac{R\theta}{v}\left
(\frac{1}{\Theta}-\frac{1}{\theta}\right)\Theta_x\frac{\bar
v_x}{\bar v}+F\frac{\bar v_x}{\bar v}.\label{3.22}
\end{eqnarray}
Because $v|_{x=0}=V|_{x=0}=v_-$, we can get
$$
\left (\frac{\bar{v}_x}{\bar{v}}\right)^2(0,t)=\left
(\frac{v_x}{v}-\frac{V_x}{V}\right)^2(0,t)=\frac{1}{s^2}\left
(\frac{u_x}{v_-}-\frac{U_x}{v_-}\right)^2(0,t)=\frac{\psi_x^2(0,t)}{s^2v_-^2}.
$$
Use Cauchy-Schwarz inequality to get
\begin{equation}
\int_0^t\left (\frac{\bar{v}_x}{\bar{v}}\right)^2(0,\tau)d\tau\leq
C\int_0^t\psi_x^2(0,\tau)d\tau\leq
C\int_0^t\left(\epsilon^{-1}\|\psi_x\|^2+\epsilon\|\psi_{xx}\|^2\right)d\tau.\label{3.23}
\end{equation}
On the other hand if we integrate (\ref{3.22}) in $R_+\times(0,t)$,
 (\ref{3.22}) is changed to
\begin{eqnarray*}
&&\int_{\R_+}\left(\frac{\mu}{2}\left (\frac{\bar v_x}{\bar
v}\right)^2-\psi\frac{\bar v_x}{\bar
v}\right)dx-\int_{\R_+}\left(\frac{\mu}{2}\left (\frac{\bar
v_{x}(x,0)}{\bar v(x,0)}\right)^2-\psi_0\frac{\bar v_{x}(x,0)}{\bar v(x,0)}\right)dx\\
&&\quad+\int_0^t\int_{\R_+}\left(\frac{R\theta}{v}\left (\frac{\bar
v_x}{\bar v}\right)^2+\left(\psi\frac{\bar v_t}{\bar
v}\right)_x-s\mu\left
(\frac{\bar v_x}{\bar v}\right)_x\frac{\bar v_x}{\bar v}\right)dxd\tau\nonumber\\
&&\leq C\epsilon^{-1}\left(\int_0^t\ \|(\zeta_x,\psi_x)\|^2\
d\tau+\int_0^t\ \int_0^{+\infty}\ \Theta_x^2(\varphi^2+\zeta^2)\ dx\
d\tau\right)\\&&+C\epsilon^{-1}\int_0^t\int_0^{+\infty}U_x^2\varphi^2
dxd\tau+C\epsilon^{-1}\int_0^t\int_0^{+\infty}|F|^2dxd\tau+\epsilon\int_0^t(\|\frac{\bar{v}_x}{\bar{v}}\|^2+\|\psi_{xx}\|^2)\
d\tau.\end{eqnarray*} Furthermore (\ref{3.22}) can be change to the
following inequality
\begin{eqnarray}\label{3.24}&&\int_0^t\|\frac{\bar{v}_x}{\bar{v}}\|^2\
d\tau+\|\frac{\bar{v}_x}{\bar{v}}\|^2-C\epsilon^{-1}\|\psi\|^2-C\|\psi_0\|^2-C\int_0^{+\infty}\frac{\bar{v}_x}{\bar{v}}(x,0)^2\
dx\nonumber\\ &&\leq C\epsilon^{-1}\left(\int_0^t\
\|(\zeta_x,\psi_x)\|^2\ d\tau+\int_0^t\ \int_0^{+\infty}\
\Theta_x^2(\varphi^2+\zeta^2)\ dx\
d\tau+C(\delta_0)\right)\nonumber\\
&&\quad+\epsilon\int_0^t\left(\|\frac{\bar{v}_x}{\bar{v}}\|^2+\|\psi_{xx}\|^2\right)\
d\tau+C\|\varphi_{0x}\|^2.
\end{eqnarray}
Because $C_1(\varphi_x^2)-C_2V_x^2\leq
(\frac{\bar{v}_x}{\bar{v}})^2\leq C_3\varphi_x^2+C_4V_x^2$
($C_1,C_2,C_3,C_4$ stands for constants about $v$), combine with
Lemma \ref{lem3.1}--\ref{lem3.2} we can find a small $\epsilon $
such that we change (\ref{3.24}) to
\begin{eqnarray}\label{3.25}\int_0^t\|\varphi_x\|^2\
d\tau+\|\varphi_x\|^2&\leq&
C\|\varphi_{0x}\|^2+C\epsilon^{-1}\|(\varphi_0,\psi_0,\zeta_0)\|^2
+\int_0^t\epsilon\|\psi_{xx}\|^2d\tau\nonumber\\
&&\quad+\int_0^tC\epsilon^{-1}\|(\psi_x,\zeta_x)\|^2d\tau+C(\delta_0).
\end{eqnarray}
 So we finish this lemma.$\Box$

 From Lemma \ref{lem3.1} to Lemma \ref{lem3.4} we know when $\delta_0$ and $\|(\varphi_0,\psi_0,\zeta_0)\|$ suitably small there exist a suitably
 small positive
 constant $\delta$
such that
$$\|(\varphi,\psi,\zeta)\|^2+\int_0^t\|(\psi_x,\zeta_x)\|^2d\tau\leq
 C\delta,$$ and $$\|(\varphi_x,\psi_x,\zeta_x)\|^2+\int_0^t\|(\psi_{xx},\zeta_{xx})\|^2\leq
 C.$$ Then we can get $C_5\leq |v|\leq
 C_6$ and $C_7\leq |\theta|\leq C_8$ when $\delta$ small, here $C_5$, $C_6,$ $C_7$ and
 $C_8$ are constants independent of $v$ and $\theta$ .
    When  combine with Lemma \ref{lem3.1}--\ref{lem3.4} we can get
 (\ref{2.1}) in
 Proposition 2.2 .

To finish Theorem \ref{thm1.1} now we will proof $\sup_{x\in
\R_+}|(\varphi,\psi,\zeta)|\to 0,\ as\ t\to \infty.$

Because $\int_0^{+\infty}\partial_x$(\ref{1.12})$_1\times
2\varphi_x\ dx$
 equals to
\begin{equation}
s\varphi_x^2(0,t)=2\int_0^\infty\varphi_x\psi_{xx}dx-\frac{d}{dt}\|\varphi_x\|^2,\label{3.26}
\end{equation}
use  Cauchy-Schwarz inequality we get
$$2\int_0^\infty\varphi_x\psi_{xx}dx \leq C\left(\|\varphi_x\|^2+\|\psi_{xx}\|^2\right),$$
according to Lemma \ref{lem3.3}--\ref{lem3.4} and (\ref{3.23}) to
get
\begin{equation}
\int_0^\infty\varphi_x^2(0,t)dt\leq
C\epsilon^{-1}\left(C(\delta_0)+\|(\varphi_0,\psi_0,\zeta_0)\|_1^2\right)+C\epsilon,\label{3.27}
\end{equation}
again using Lemma \ref{lem3.3}--\ref{lem3.4} and (\ref{3.26}), then
from (\ref{3.27}) we get
\begin{eqnarray}
&&\int_0^\infty\left|\frac{d}{dt}\|\varphi_x(t)\|^2\right|dt\nonumber\\
&&\leq
C\int_0^\infty\varphi_x^2(0,t)dt+C\int_0^\infty\left(\|\varphi_x\|^2+\|\psi_{xx}\|^2\right)dt\nonumber\\
&&\leq
C\epsilon^{-1}\left(C(\delta_0)+\|(\varphi_0,\psi_0,\zeta_0)\|_1^2\right)+C\epsilon.\label{3.28}
\end{eqnarray}
Similar as above, from Lemma \ref{lem3.1}$-$\ref{lem3.4} and combine
with Sobolev inequality we get
\begin{equation}
\int_0^\infty\left(\left|\frac{d}{dt}\|\psi_x(t)\|^2\right|+\left|\frac{d}{dt}\|\zeta_x(t)\|^2\right|\right)d\tau\leq
C\epsilon^{-1}_4\left(C(\delta_0)+\|(\varphi_0,\psi_0,\zeta_0)\|_1^2\right)+C\epsilon.\label{3.29}
\end{equation}

It means
$$
\|(\varphi,\psi,\zeta)(t)\|_{L^\infty}^2\leq
2\|(\varphi,\psi,\zeta)(t)\|\|(\varphi_x,\psi_x,\zeta_x)(t)\|\to0\quad\mbox{when}\quad
t\to\infty.
$$
So we finish Theorem \ref{thm1.1}.$\Box$

\end{document}